\documentclass[a4paper,11pt]{article}
\usepackage[latin1]{inputenc}
\usepackage[francais,english]{babel}    
\usepackage{amssymb}
\usepackage{amsmath}
\usepackage{amsthm}
\usepackage{multicol}
\usepackage{graphicx}
\usepackage{color}
\usepackage[T1]{fontenc}
\usepackage{yfonts}
\usepackage{placeins}
\usepackage{array}
\usepackage{dsfont}
\usepackage{stmaryrd}
\usepackage{verbatim}
\usepackage{caption}
\usepackage{transparent}
\usepackage{bbold}
\usepackage{hyperref}
\usepackage{enumitem}
\usepackage{tikz}
\usetikzlibrary{patterns}
\usepackage{tikz-3dplot}

\newtheorem{thm}{Theorem}[section]
\newcommand{\bthm}{\begin{thm}}
\newcommand{\ethm}{\end{thm}}

\newtheorem{thmi}{Theorem}
\newcommand{\bthmi}{\begin{thmi}}
\newcommand{\ethmi}{\end{thmi}}

\newtheorem{cori}[thmi]{Corollary}
\newcommand{\bcori}{\begin{cori}}
\newcommand{\ecori}{\end{cori}}

\newtheorem{mthm}{Theorem}
\newcommand{\bmthm}{\begin{mthm}}
\newcommand{\emthm}{\end{mthm}}
\newtheorem{mcor}[mthm]{Corollary}
\newcommand{\bmcor}{\begin{mcor}}
\newcommand{\emcor}{\end{mcor}}

\newtheorem*{conj}{Conjecture}
\newcommand{\bconj}{\begin{conj}}
\newcommand{\econj}{\end{conj}}

\newtheorem*{question}{Question}
\newcommand{\bq}{\begin{question}}
\newcommand{\eq}{\end{question}}

\newtheorem*{thn}{Theorem}
\newcommand{\bthn}{\begin{thn}}
\newcommand{\ethn}{\end{thn}}

\newtheorem{exo}{Exercise}
\newcommand{\bex}{\begin{exo}}
\newcommand{\eex}{\end{exo}}

\newtheorem{sol}{Solution}
\newcommand{\bsol}{\begin{sol}}
\newcommand{\esol}{\end{sol}}

\newtheorem{pro}[thm]{Proposition}
\newcommand{\bpro}{\begin{pro}}
\newcommand{\epro}{\end{pro}}

\newtheorem{cor}[thm]{Corollary}
\newcommand{\bcor}{\begin{cor}}
\newcommand{\ecor}{\end{cor}}

\newtheorem{lem}[thm]{Lemma}
\newcommand{\blem}{\begin{lem}}
\newcommand{\elem}{\end{lem}}

\theoremstyle{definition}

\newtheorem{defi}[thm]{Definition}
\newcommand{\bdf}{\begin{defi}}
\newcommand{\edf}{\end{defi}}

\newtheorem*{defis}{Definition}
\newcommand{\bdfs}{\begin{defis}}
\newcommand{\edfs}{\end{defis}}

\newtheorem*{rmk}{Remark}
\newcommand{\brk}{\begin{rmk} \upshape}
\newcommand{\erk}{\end{rmk}}

\newtheorem*{rmks}{Remarks}
\newcommand{\brks}{\begin{rmks} \upshape}
\newcommand{\erks}{\end{rmks}}

\newtheorem*{exe}{Example}
\newcommand{\bexe}{\begin{exe} \upshape}
\newcommand{\eexe}{\end{exe}}

\newtheorem*{exes}{Examples}
\newcommand{\bexes}{\begin{exes} \upshape}
\newcommand{\eexes}{\end{exes}}

\newtheorem*{pre}{Proof}
\newcommand{\bp}{\begin{pre} \upshape}
\newcommand{\ep}{\hfill \qed \end{pre}}
\newcommand{\epp}{\end{pre}}

\newtheorem{mainthm}{Theorem}
\newtheorem{maincor}[mainthm]{Corollary}

\newtheorem{mainconj}[mainthm]{Conjecture}

\newcommand{\beq}{\begin{eqnarray*}}
\newcommand{\eeq}{\end{eqnarray*}}

\newcommand{\beqn}{\begin{equation}}
\newcommand{\eeqn}{\end{equation}}

\newcommand{\ben}{\begin{enumerate}}
\newcommand{\een}{\end{enumerate}}
\newcommand{\bit}{\begin{itemize} \renewcommand{\labelitemi}{$\bullet$} \renewcommand{\labelitemii}{$\star$}}
\newcommand{\eit}{\end{itemize}}

\newcommand{\bfg}{
\begin{figure}[H]
\begin{center}}
\newcommand{\efg}{
\end{center}
\end{figure}
\FloatBarrier}

\newcolumntype{M}[1]{>{\raggedright}m{#1}}

\newcommand{\R}{\mathbb{R}}

\newcommand{\N}{\mathbb{N}}
\newcommand{\Z}{\mathbb{Z}}

\newcommand{\F}{\mathbb{F}}

\newcommand{\bs}{\symbol{92}}

\newcommand{\llb}{\llbracket}
\newcommand{\rrb}{\rrbracket}

\renewcommand{\dim}{\operatorname{dim}}

\newcommand{\Min}{\operatorname{Min}}

\newcommand{\bord}{\partial_\infty}

\newcommand{\Hull}{\operatorname{Hull}}
\newcommand{\rk}{\operatorname{rk}}

\newcommand{\st}{\, | \,}

\newcommand{\ra}{\rightarrow}

\newcommand{\ral}[1]{\underset{#1}{\longrightarrow}}

\newcommand{\f}{\frac}

\renewcommand{\geq}{\geqslant}
\renewcommand{\leq}{\leqslant}

\newcommand{\PSL}{\operatorname{PSL}}

\newcommand{\<}{\langle}
\renewcommand{\>}{\rangle}
\newcommand{\pif}{{+\infty}}

\newcommand{\mk}{\medskip}

\newcommand{\sign}{\begin{flushright}
Thomas Haettel \\
Institut Montpelli\'{e}rain Alexander Grothendieck \\
CNRS, Univ. Montpellier \\
thomas.haettel@umontpellier.fr
\end{flushright}}

\makeatletter
\def\Ddots{\mathinner{\mkern1mu\raise\p@
\vbox{\kern7\p@\hbox{.}}\mkern2mu
\raise4\p@\hbox{.}\mkern2mu\raise7\p@\hbox{.}\mkern1mu}}
\makeatother

\makeatletter
\def\maketitles{%
  \null
  \thispagestyle{empty}%
  \vfill
  \begin{center}\leavevmode
    \normalfont
    {\LARGE \@title\par}%
    \vskip 1.2cm
    {\large \@author\par}%
    \vskip 1.2cm
    {\large \@subtitle\par}%
    \vskip 0.8cm
    {\large \@date\par}%
  \end{center}%
  \vfill
  \null
  \cleardoublepage
  }

\def\date#1{\def\@date{#1}}
\def\author#1{\def\@author{#1}}
\def\title#1{\def\@title{#1}}
\def\subtitle#1{\def\@subtitle{#1}}
\makeatother

\newtheorem*{corG}{Corollary~\ref{cor:braid}}
\newtheorem*{thmE}{Theorem~\ref{thm:main_dagger}}

\title{Virtually cocompactly cubulated Artin-Tits groups}
\author{Thomas Haettel}
\date{\today}

\begin{document}

\selectlanguage{english}

\maketitle

\begin{center}
\begin{minipage}{0.8\textwidth}
\textsc{Abstract.} We give a conjectural classification of virtually cocompactly cubulated Artin-Tits groups (i.e. having a finite index subgroup acting geometrically on a CAT(0) cube complex), which we prove for all Artin-Tits groups of spherical type, FC type or two-dimensional type. A particular case is that for $n \geq 4$, the $n$-strand braid group is not virtually cocompactly cubulated. \end{minipage}
\end{center}

\let\thefootnote\relax

\section*{Introduction}

Groups acting geometrically on CAT(0) spaces (called CAT(0) groups), or even better on CAT(0) cube complexes (called cocompactly cubulated groups), possibly up to a finite index subgroup, enjoy a list of nice properties: they have a quadratic Dehn function, a solvable word and conjugacy problem, they have the Haagerup property, their amenable subgroups are virtually abelian and undistorted, they satisfy the Tits alternative... R.~Charney conjectures that all Artin-Tits groups are CAT(0), but very few cases are known. With D.~Kielak and P.~Schwer (see~\cite{b6}), we pursued the construction of T.~Brady and J.~McCammond (see~\cite{brady} and \cite{brady_mccammond}) to prove that for $n \leq 6$, the $n$-strand braid group is CAT(0).

In this article, we give a conjectural classification of which Artin-Tits groups are virtually cocompactly cubulated, and we prove this classification under a mild conjecture on Artin-Tits groups, which is satisfied in particular for spherical, type FC or $2$-dimensional Artin-Tits groups. Right-angled Artin groups are well-known to act cocompactly on their Salvetti CAT(0) cube complex, but there are a few more examples. This question was asked by D.~Wise for the particular case of braid groups (see~\cite[Problem~13.4]{wise_cubical_route}).

\begin{mainconj}[Classification of virtually cocompactly cubulated Artin-Tits groups] \label{conj:main}\

Let $M=(m_{ab})_{a,b \in S}$ be a finite Coxeter matrix. Then the Artin-Tits group $A(M)$ is virtually cocompactly cubulated if and only if the following two conditions are satisfied:
\ben
\item for each pairwise distinct $a,b,c \in S$ such that $m_{ab}$ is odd, either $m_{ac}=m_{bc}=\infty$ or $m_{ac}=m_{bc}=2$, and
\item for each distinct $a,b \in S$ such that $m_{ab}$ is even and different from $2$, there is an ordering of $\{a,b\}$ (say $a < b$) such that, for every $c \in S \bs \{a,b\}$, one of the following holds:
\bit
\item $m_{ac}=m_{bc}=2$,
\item $m_{ac}=2$ and $m_{bc}=\infty$,
\item $m_{ac}=m_{bc}=\infty$, or
\item $m_{ac}$ is even and different from $2$, $a<c$ in the ordering of $\{a,c\}$, and $m_{bc}=\infty$.
\eit
\een
\end{mainconj}

In particular, typical examples of cocompactly cubulated Artin-Tits groups are the following types:
\bit
\item right-angled Artin groups, i.e. such that $\forall a,b \in S, m_{ab} \in \{2,\infty\}$,
\item dihedral Artin groups, i.e. such that $|S|=2$,
\item ``even stars'' Artin groups, i.e. such that there exists a ``central vertex'' $a_0 \in S$ such that $\forall a,b \in S \bs \{a_0\}, m_{ab}=\infty$ and $\forall a \in S\bs \{a_0\}, m_{aa_0}$ is even.
\eit

\begin{figure}[!h]
\def\svgwidth{4cm}
\center
\input{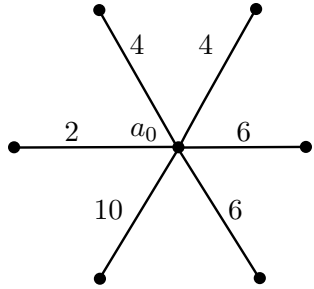}
\caption{An example of the Coxeter graph of an ``even star'' cocompactly cubulated Artin-Tits group, with central vertex $a_0$}
\label{fig:examples}
\end{figure}

You can see in Figure~\ref{fig:examples} an example of the Coxeter graph of an even star Artin-Tits group. In that figure, all the edges labeled $\infty$ are not drawn.

The most general picture of an arbitrary cocompactly cubulated Artin-Tits group comes roughly from combining dihedral Artin groups and even stars Artin groups, in a right-angled-like fashion.

\begin{figure}[!h]
\def\svgwidth{11cm}
\center
\input{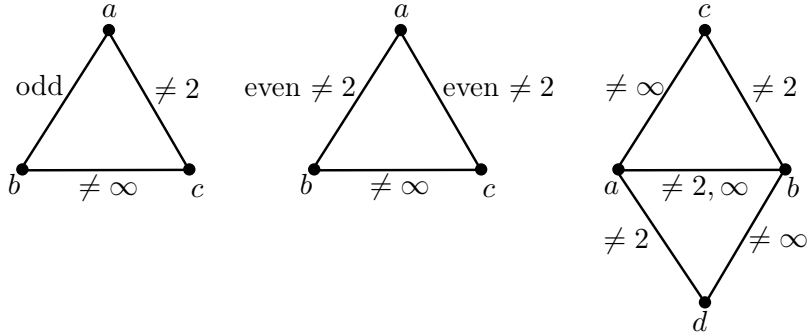}
\caption{Local obstructions in Conjecture~\ref{conj:main2}}
\label{fig:obstructions}
\end{figure}

Another way to state Conjecture~\ref{conj:main} is by describing local obstructions in the Coxeter matrix $M$, see also Figure~\ref{fig:obstructions}. In particular, standard parabolic subgroups of rank $3$ and $4$ should determine if an Artin-Tits group is virtually cocompactly cubulated or not. Recall that the rank of an Artin group (or a standard parabolic subgroup) is the number of standard generators.

\begin{mainconj}[Reformulation of Conjecture~\ref{conj:main}] \label{conj:main2}
Let $M=(m_{ab})_{a,b \in S}$ be a finite Coxeter matrix. Then the Artin-Tits group $A(M)$ is not virtually cocompactly cubulated if and only one of the following occurs:
\bit
\item there exist $3$ pairwise distinct $a,b,c \in S$ such that $m_{ab}$ is odd, $m_{ac} \neq \infty$ and $m_{bc} \neq 2$,
\item there exist $3$ pairwise distinct $a,b,c \in S$ such that $m_{ab}$ and $m_{ac}$ are even numbers different from $2$, and $m_{bc} \neq \infty$, or
\item there exist $4$ pairwise distinct $a,b,c,d \in S$ such that $m_{ab} \not\in \{2,\infty\}$, $m_{ac}, m_{bd} \neq \infty$ and $m_{ad},m_{bc} \neq 2$.
\eit
\end{mainconj}

One implication of Conjecture~\ref{conj:main}, namely the cubulation of Artin-Tits groups satisfying the two conditions, is proven in full generality in this article.

\begin{mainthm} \label{thm:main1}
Let $M=(m_{ab})_{a,b \in S}$ be a finite Coxeter matrix satisfying the two conditions of Conjecture~\ref{conj:main}. Then $A(M)$ is cocompactly cubulated.
\end{mainthm}

The converse implication of Conjecture~\ref{conj:main}, namely to show that the two conditions are necessary to be virtually cocompactly cubulated, is proven under the following mild assumption. Let $M=(m_{ab})_{a,b \in S}$ be a finite Coxeter matrix. We say that the Artin-Tits group $A(M)$ satisfies property $(\dagger)$ if
$$\forall s \in S, \forall n \geq 1, Z_{A(M)}(s^n) = Z_{A(M)}(s).$$

\begin{mainthm} \label{thm:main2}
Let $M=(m_{ab})_{a,b \in S}$ be a finite Coxeter matrix such that $A(M)$ satisfies property $(\dagger)$. If $A(M)$ is virtually cocompactly cubulated, then $A(M)$ satisfies the two conditions of Conjecture~\ref{conj:main}.
\end{mainthm}

It is conjectured that all Artin-Tits groups satisfy property $(\dagger)$, it is notably a very restricted consequence of Property $(\star)$ in~\cite{godelle_cat0}. In particular, it is true as soon as the Deligne complex can be endowed with a piecewise Euclidean CAT(0) metric. This condition is therefore true for Artin-Tits groups of type FC (i.e. every complete subgraph spans a spherical subgroup), with the cubical metric on the Deligne complex (see~\cite[Theorem~A]{charney_davis_kpi1}). It is also true if the Artin-Tits group $A(M)$ is such that any irreducible spherical parabolic subgroup has rank at most $2$ (which is slightly more general than $2$-dimensional), in which case the Moussong metric on the Deligne complex is CAT(0) (see~\cite[Theorem~A]{charney_davis_kpi1}). Note that the Moussong metric on the Deligne complex is conjectured to be CAT(0) for all Artin-Tits groups (see~\cite[Conjecture~3]{charney_davis_kpi1}).

\begin{mainthm} \label{thm:main_dagger}
Conjecture~\ref{conj:main} holds for any Artin-Tits group satisfying property $(\dagger)$. In particular, Conjecture~\ref{conj:main} holds for Artin-Tits groups of type FC, and for Artin-Tits groups whose irreducible spherical parabolic subgroups have rank at most $2$.
\end{mainthm}

One should note that the two conditions in Conjecture~\ref{conj:main} imply that the Artin-Tits group is of type FC. In particular, the following is a consequence of Conjecture~\ref{conj:main}.

\begin{mainconj}[Consequence of Conjecture~\ref{conj:main}] \label{conj:mainFC}
If an Artin-Tits group is virtually cocompactly cubulated, then it is of type FC.
\end{mainconj}

In the particular case or Artin groups of spherical type, the condition is much simpler.

\begin{maincor}[Classification of virtually cocompactly cubulated Artin-Tits groups of spherical type] \label{cor:spherical}
Let $A$ be an Artin-Tits group of spherical type. Then $A$ is virtually cocompactly cubulated if and only if every irreducible parabolic subgroup of $A$ has rank at most $2$.
\end{maincor}

In particular, this gives a very simple answer for braid groups.

\begin{maincor}[Cubulation of braid groups] \label{cor:braid}
The $n$-strand braid group $B_n$, or its central quotient $B_n/Z(B_n)$, is virtually cocompactly cubulated if and only if $n \leq 3$.
\end{maincor}

However, according to B.~Bowditch (see~\cite{bowditch_coarse_median}), all mapping class groups, including braid groups, are coarse median, which implies that their asymptotic cones \og look like \fg\ asymptotic cones of CAT(0) cube complexes: they are not cocompactly cubulated, but \og look cubical \fg\ on a large scale.

\mk

Concerning proper actions of Artin groups on CAT(0) cube complexes, even the following question is still open. 

\bq[Charney~\cite{charney_problems}, Wise~\cite{wise_cubical_route}]
Does the $4$-strand braid group $B_4$ have a metrically proper action on a CAT(0) cube complex ?
\eq

During the proof, we also prove the following cubulation results, of independent interest. See Theorem~\ref{thm:cubulation of multiple centralizer} and Proposition~\ref{pro:cubulation_central_quotient} fore more precise versions.

\begin{mainthm}[Cubulation of normalizers and centralizers] \label{thm:stability_normalizers} Let $G$ be a cocompactly cubulated group, and let $A$ be an abelian subgroup of $G$. Then $A$ has a finite index subgroup $A_0$ such that $N_G(A_0)$ is cocompactly cubulated, and $Z_G(A_0)$ has finite index in $N_G(A_0)$. \end{mainthm}

\begin{mainthm}[Cubulation of central quotients] \label{thm:stability_central_quotients} Let $G$ be a cocompactly cubulated group, and let $A$ be a central, convex-cocompact subgroup of $G$. Then $G/A$ is cocompactly cubulated. \end{mainthm}

In an earlier version of this article, we obtained only Theorem~\ref{thm:main2} for cocompactly cubulated Artin groups, without the virtual part. J.~Huang, K.~Jankiewicz and P.~Przytycki, simultaneously to this earlier version and independently, proved Theorem~\ref{thm:main2} for $2$-dimensional Artin groups with the virtual part (see~\cite{huang_jankiewicz_przytycki}). In particular, they showed that a $2$-dimensional Artin group is cocompactly cubulated if and only if it is virtually cocompactly cubulated. 

\mk

Concerning Coxeter groups, Niblo and Reeves proved (see~\cite{niblo_reeves}) that every Coxeter group acts properly on a locally finite CAT(0) cube complex. Caprace and M\"uhlherr proved (see~\cite{caprace_muhlherr}) that this action is cocompact if and only if the Coxeter diagram does not contain an affine subdiagram of rank at least $3$.

\mk

O.~Varghese recently described (see~\cite{varghese}) a group-theoretic condition ensuring that any (strongly simplicial) isometric action on a CAT(0) cube complex has a global fixed point. This condition is notably satisfied by $Aut(\F_n)$ for $n \geq 1$.

\medskip

\textbf{Outline of the proof} The rough idea is to study the CAT(0) visual angle between maximal abelian subgroups in Artin groups. Using a result from J.~Crisp and L.~Paoluzzi (see~\cite{crisp_paoluzzi}), we show that if $a,b$ are the standard generators of the $3$-strand braid group acting on some CAT(0) space, then the translation axes for $a$ and $ababab$ form an acute visual angle at infinity.

\mk

On the other hand, we show that the translation axes of elements in maximal abelian subgroups of a group acting geometrically on a CAT(0) cube complex, with finite intersection, form an obtuse visual angle at infinity. This is the source of the non-cubicality results. This uses a flat torus theorem for maximal abelian subgroups of cocompactly cubulated groups by Wise and Woodhouse (see~\cite{wise_woodhouse} and Theorem~\ref{thm:flat_torus}).

\medskip

\textbf{Acknowledgments:} The author would like to thank very warmly Jingyin Huang, Kasia Jankiewicz and Piotr Przytycki for their precious help to find and correct several mistakes in previous versions of this article. The author would also like to thank Daniel Wise for inspiring discussions. The author would also like to thank Eddy Godelle and Luis Paris, for very interesting discussions on Artin groups. The author would also like to thank Nir Lazarovich and Anthony Genevois for noticing several typos and a few mistakes. Finally, the author would like to thank several anonymous referees for pointing out mistakes in previous versions of this article, and for many useful comments that helped in particular to improve the exposition.
\section{Definitions and notations}

\subsection{Artin groups}

For $p \in \N$, let $w_p$ denote the word $w_p(a,b)=aba...ba$ of length $p$. Let $S$ be a finite set, and let $\Gamma$ be a graph with vertex set $S$ and edges labeled in $\N_{\geq 2}$. The \emph{Artin-Tits group} $A(\Gamma)$ is defined by the following presentation:
$$ A(\Gamma) = \left\< s \in S \st w_p(s,t) = w_p(t,s) \mbox{ for each edge $\{s,t\}$ labeled $p$}\right\>.$$

If $S=\{a,b\}$, then $A(\Gamma)$ is called a dihedral Artin group, and we will denote it by $A(p)$, where $p$ is the label of the edge $\{a,b\}$ (or $p=\infty$ if there is no edge). For instance, $A(2) \simeq \Z^2$ and $A(\infty) \simeq \F_2$.

If $a$ and $b$ are different elements of $S$, then the subgroup of $A(\Gamma)$ spanned by $a$ and $b$ is isomorphic to the dihedral group $A(p)$, where $p$ is the label of the edge $\{a,b\}$. If $p \not\in \{2,\infty\}$, the center of $A(p)$ is the infinite cyclic group spanned by $z_{ab}=w_q(a,b)$, where $q=2p$ if $p$ is odd, and $q=p$ if $p$ is even (see~\cite{brieskorn_saito} and \cite{deligne} for the center of spherical Artin groups). If $p=2$, the dihedral group $A(2)$ is free abelian, we will denote $z_{ab}=ab$ in this case.

\subsection{CAT(0) cube complexes}

A finite dimensional cube complex $X$ is naturally endowed with two natural distances, defined piecewise on cubes: the $L^1$ distance $d_1$ and the $L^2$ distance $d_2$ (each edge has length $1$). Throughout the paper, unless we want to use both distances, we will mainly use the $L^1$ distance $d_1$ and will simply denote it $d$.

A cube complex $X$ is called CAT(0) if the $d_2$ distance is CAT(0), or equivalently if the $d_1$ distance is median (see section~\ref{subsec:median}). A discrete group $G$ is called cocompactly cubulated if it acts geometrically, i.e. properly and cocompactly by cubical isometries, on a CAT(0) cube complex.

Let us recall the fundamental local-to-global property for CAT(0) spaces.

\bthm[Cartan-Hadamard] A metric space is CAT(0) if and only if it is simply connected and locally CAT(0).
\ethm

Let us recall Gromov's combinatorial criterion to show that a cube complex is locally CAT(0).

\bthm[Gromov, see~\cite{gromov_hyperbolic_groups}] \label{thm:gromov criterion}
A cube complex $X$ is locally CAT(0) if and only if, for any $3$ cubes $Q$,$Q'$,$Q''$ of $X$, which pairwise intersect in codimension $1$ and intersect globally in codimension $2$, they are codimension $1$ faces of some cube of $X$.
\ethm

In a CAT(0) cube complex $X$, a \emph{hyperplane} $H$ denotes the orthogonal (with respect to the CAT(0) metric $d_2$) of some edge $[x,y]$ at its midpoint (i.e. the set of points of $X$ whose projection on $[x,y]$ equals the midpoint), we denote it $H=[x,y]^\perp$ (the hyperplane $H$ can also be described as a union of midcubes, or as an equivalence class of edges, see~\cite{sageev}). Each hyperplane divides $X$ into two connected components, the closures of which are called \emph{half-spaces} and denoted by $H^+$ and $H^-$. An automorphism $g$ of $X$ is said to \emph{skewer} the half-space $H^+$ if $g \cdot H^+ \subsetneq H^+$. By skewering $H$, we mean skewering $H^+$ or $H^-$.

We say that two hyperplanes $H,H'$ \emph{cross} if $H^+$ and $H^-$ intersect $H'^+$ and $H'^-$.

If $x,y$ are vertices of a CAT(0) cube complex $X$, then $d_1(x,y)$, also called the \emph{combinatorial distance} between $x$ and $y$, coincides with the number of hyperplanes separating $x$ and $y$. An automorphism $g$ of $X$ is called \emph{combinatorially hyperbolic} if $g$ preserves a combinatorial ($d_1$) geodesic, on which it acts by a nontrivial translation. An action of a group $G$ by cubical automorphisms on $X$ is called \emph{combinatorially semisimple} if every element of $G$ is either combinatorially hyperbolic or fixes a vertex.

If $X$ is a cube complex, we can divide naturally each $d$-cube into $2^d$ smaller cubes, getting a new cube complex (up to rescaling the metric by $2$) called the \emph{cubical subdivision} of $X$.

\bthm[Haglund, see~\cite{haglund}] \label{thm:haglund} Let $G$ be a group acting by cubical automorphisms on a CAT(0) cube complex $X$. Then $G$ acts combinatorially semisimply on the cubical subdivision of $X$. \ethm

In particular, if $G$ acts properly on $X$ and $g \in G$ has infinite order, then $g$ acts as a combinatorial hyperbolic isometry of the cubical subdivision of $X$.

\medskip

If $g$ is a cubical isometry of a CAT(0) cube complex, its (combinatorial) \emph{translation length} is
$$\delta_g = \min_{x \in X^{(0)}} d_1(x,g \cdot x).$$

If $A$ is a subgroup of $G$, its (combinatorial) \emph{minimal set} is
$$\Min(A) = \{x \in X \st \forall a \in A, d_1(x,a \cdot x) = \delta_a\}.$$
If $g \in G$, we will simply denote $\Min(g)$ instead of $\Min(\<g\>)$.

\brk Note that, according to Theorem~\ref{thm:haglund}, for any cubical automorphism $g$ of a CAT(0) cube complex $X$, we have $\Min(g) \neq \emptyset$. Furthermore, up to passing to the cubical subdivision of $X$, we have $\Min^{(0)}(g) = \Min(g) \cap X^{(0)} \neq \emptyset$. \erk

\brk Also note that $\Min(g)$ need not be convex for the $d_1$ distance, nor need it be a cube subcomplex: consider for instance $X=\R^2$, with the standard Cayley square complex structure of $\Z^2$, and let $g:(x,y) \mapsto (y+1,x+1)$. Then $\delta_g=2$ and $\Min(g)=\{(x,y) \in \R^2 \st |x-y| \leq 1\}$ is not a cube subcomplex. \erk

We can also give a slight variation of this result, by considering an appropriate power of $g$ instead of considering the cubical subdivision.

\bpro \label{pro:min_vertex}
Let $X$ be a CAT(0) cube complex of dimension at most $D$, and let $g$ denote a cubical automorphism of $X$. Then $\Min(g^{2D!}) \cap X^{(0)} \neq \emptyset$.
\epro

\bp
Let $X'$ denote the cubical subidivision of $X$. According to Theorem~\ref{thm:haglund}, either $g$ fixes a vertex of $X'$ or $g$ is a combinatorially hyperbolic isometry of $X'$.

\mk

Assume that $g$ fixes a vertex of $X'$, i.e. $g$ stabilizes a cube $Q$ of $X$. If $Q$ is a $k$-cube, the isometry group of $Q$ is isomorphic to $\frak{S}_k \ltimes \left(\Z/2\Z\right)^k$. Each isometry of $Q$ has an order which divides $2k!$, and so it also divides $2D!$. Hence $g^{2D!}$ fixes each vertex of $Q \cap X^{(0)}$.

\mk

Assume that $g$ is a combinatorially hyperbolic isometry of $X'$. So there exists a vertex $v$ of $X'$ such that $\{g^n \cdot v, n \in \Z\}$ lie on a combinatorial geodesic in $X'$. The vertex $v$ of $X'$ corresponds to the midpoint of a cube $Q$ of $X$. Choose $v$ such that the dimension $k$ of $Q$ is minimal. By minimality of $k$, each hyperplane $H$ of $X$ (seen as a union of midcubes) containing $v$ also contains $g \cdot v$. So $g$ preserves the set of $k$ hyperplanes of $X$ containing $v$. As a consequence, the element $g^{D!}$ fixes each hyperplane of $X$ containing $v$, and $g^{2D!}$ furthermore preserves the orientations of each of these $k$ hyperplanes. This implies that for each vertex $x$ of $Q \cap X^{(0)}$, we have $d(x,g^{2D!} \cdot x)=d(v,g^{2D!} \cdot v)$. Hence $x \in \Min(g^{2D!})$, so $Q \cap X^{(0)} \subset \Min(g^{2D!})$.
\ep

If $X$ is a proper CAT(0) cube complex, we will denote by $\bord X$ its visual (CAT(0)) boundary at infinity: it is endowed with the visual distance $\sphericalangle$. Each combinatorially hyperbolic isometry $g$ of $X$ is CAT(0) hyperbolic, and has a unique attracting fixed point $g(+\infty) \in \bord X$. See~\cite{bridson_haefliger} for details on general CAT(0) spaces and isometries.

\subsection{Median algebras} \label{subsec:median}

A \emph{median algebra} is a set $M$ endowed with a symmetric map $\mu:M^3 \ra M$, called the \emph{median}, satisfying the following
\begin{align*}
&\forall a,b \in M, \mu(a,a,b)=a \\
&\forall a,b,c,d,e \in M, \mu(a,b,\mu(c,d,e))=\mu(\mu(a,b,c),\mu(a,b,d),e). \addtocounter{equation}{1}\tag{\theequation} \label{eqn:median}
\end{align*}

In a metric space $M$, the interval between $a \in M$ and $b \in M$ denotes $I(a,b)=\{c \in M \st d(a,c)+d(c,b)=d(a,b)\}$. A metric space $M$ is called \emph{metric median} if
$$ \forall a,b,c \in M, I(a,b) \cap I(b,c) \cap I(a,c) =\{\mu(a,b,c)\},$$
which implies that $\mu$ is a median (see~\cite{bandelt_hedlikova}).

Median algebras and CAT(0) cube complexes are highly related, as proved by Chepoi.

\bthm[\cite{chepoi_median_cat0}] \label{thm:chepoi}
A connected graph, endowed with its combinatorial distance, is metric median if and only if it is the $1$-skeleton of a CAT(0) cube complex.
\ethm

Recall that the \emph{rank} of a subset $I$ of a median algebra $M$ is the supremum of all $n \in \N$ such that the $n$-cube $\{0,1\}^n$ has a median-preserving embedding into $I$. Starting with a more general median algebra, one has the following. 

\bthm[\cite{nica} and \cite{chatterji_niblo}] \label{thm:CCC from median}
Let $M$ be a median algebra with intervals of rank at most $D$, and let $G$ be a group of automorphisms of $M$. There exists a CAT(0) cube complex $X(M)$, with vertex set $X(M)^{(0)}=M$, of dimension at most $D$, on which $G$ acts as a group of cubical automorphisms. \ethm

If $X$ is a CAT(0) cube complex and $x,y \in X$, let $I(x,y)=\{z \in X \st d_1(x,z) + d_1(z,y)=d_1(x,y)\}$ denote the $d_1$ interval between $x$ and $y$. A subset $Y \subset X$ is said to be \emph{convex} if for every $x,y \in Y$, we have $I(x,y) \subset Y$. If $Y \subset X$, the \emph{$d_1$ convex hull} of $Y$ is the smallest convex subset of $X$ containing $Y$, denoted $\Hull(Y)$.

\mk

The median $\mu :X^3 \ra X$ is defined by
$$ \forall x,y,z \in X, I(x,y) \cap I(y,z) \cap I(z,x) = \{\mu(x,y,z)\}.$$

It is well-known that $\mu$ is $1$-Lipschitz for $d_1$ with respect to each of the three variables. For the convenience of the reader, we give a short proof when restricted to the set of vertices.

\blem \label{lem:lipschitz}
The map $\mu : {X^{(0)}}^3 \ra X^{(0)}$ is $1$-Lipschitz with respect to the three variables, for the distance $d_1$.
\elem

\bp
Let $x,x' \in X^{(0)}$ be adjacent vertices of $X$ and $y,z \in X^{(0)}$. Let $m=\mu(x,y,z)$ and $m'=\mu(x',y,z)$. Assume that $H$ is a hyperplane of $X$ separating $m$ and $m'$, say $m \in H^+$ and $m' \in H^-$. Since $m,m' \in I(y,z)$ we deduce that $x \in H^+$ and $x' \in H^-$. So there exists at most $1$ hyperplane separating $m$ and $m'$: $d_1(m,m') \leq 1$. This implies that $\mu$ is $1$-Lipschitz. 
\ep

If $I \subset \R$ is an interval, a map $c:I \ra X$ is called monotone if
$$\forall s \leq t \leq u \in I, \mu(c(s),c(t),c(u)) = c(t).$$

If $C \subset X$ is a non-empty convex subset, there exists a unique map $\pi_C : X \ra C$, called the gate projection onto $C$ (see for instance~\cite{hhg1}), such that
$$ \forall x \in X, \forall c \in C, \mu(x,\pi_C(x),c)=\pi_C(x).$$

We now state a Lemma which will be used later.

\blem \label{lem:median inter convex}
Let $X$ denote a CAT(0) cube complex, and consider the median algebra $(X^{(0)},\mu)$. Let $M \subset X^{(0)}$ denote a median subalgebra, and let $C \subset X$ denote a convex subset. Then $M \cap C$ is a convex subset of $M$.\elem

\bp
Fix $x,y \in M \cap C$, we want to prove that the interval $I_M(x,y)$ between $x$ and $y$ in the median algebra $M$ is contained in $M \cap C$. Let $I_X(x,y)$ denote the interval between $x$ and $y$ in $X$, we have $I_M(x,y) = I_X(x,y) \cap M$. Since $C$ is convex in $X$, we have $I_X(x,y) \subset C$. Hence $I_M(x,y) \subset C \cap M$. So $M \cap C$ is a convex subset of $M$.
\ep

If a group $G$ acts by cubical isometries on a CAT(0) cube complex $X$, the action is said to be \emph{median minimal} if $X^{(0)}$ is the smallest $G$-invariant non-empty median subalgebra of $X^{(0)}$.

\section{Cubulation of centralizers}

In this section, we prove the following result on cubulation of centralizers, which is more precise than Theorem~\ref{thm:stability_normalizers} stated in the introduction.

\begin{thm} \label{thm:cubulation of multiple centralizer} Let $G$ be a group acting geometrically by isometries on a locally finite CAT(0) cube complex $X$ of dimension at most $D$. Let $A$ be an abelian subgroup of $G$ such that every element of $A$ is the $2D!^\text{th}$ power of a combinatorially hyperbolic isometry in $G$. Then the normalizer $N_G(A)$ of $A$ acts geometrically on the locally finite CAT(0) cube complex $X(\Min(A)^{(0)})$ of dimension at most $D$ associated to the median subalgebra $\Min(A)^{(0)}$. Furthermore, the centralizer $Z_G(A)$ has finite index in $N_G(A)$.
\end{thm}

\brk It is not always true that $Z_G(g)$ acts cocompactly on $\Min(g^{2D!})^{(0)}$: consider for instance $X=\R^2$, with the standard Cayley square complex structure of $\Z^2$, and let $g:(x,y) \mapsto (y+1,x+1)$ and $h:(x,y) \mapsto (x+1,y)$. Consider the group $G$ spanned by $g$ and $h$. We have $Z_G(g)=\<g\>\simeq \Z$, but $\Min(g^4)=\R^2$. \erk

\blem \label{lem:case M(g)} Let $X$ be a CAT(0) cube complex of dimension at most $D$, and let $g$ be a combinatorially hyperbolic isometry of $X$. Then for any $x \in \Min(g)^{(0)}$, and for any hyperplane $H$ separating $x$ and $g^{D!} \cdot x$, $g^{D!}$ skewers $H$. \elem

\bp
Fix $x \in \Min(g)^{(0)}$, and let $H$ be a hyperplane such that $x \in H^-$ and $g^{D!} \cdot x \in H^+$. Note that, since there is a combinatorial geodesic from $x$ to $g^{2D!} \cdot x$ through $g^{D!} \cdot x$ crossing $H$ and $g^{D!} \cdot H$, we know that $g^{D!} \cdot H \neq H$.

Assume that for every $0 \leq i < j \leq D$ we have $g^i \cdot H \cap g^j \cdot H \neq \emptyset$. Since $g^{D!} \cdot H \neq H$, for every $0 \leq i < j \leq D$, we have $g^i \cdot H \neq g^j \cdot H$ and $g^i \cdot H \cap g^j \cdot H \neq \emptyset$ so $g^i \cdot H$ and $g^j \cdot H$ cross. Hence the $D+1$ hyperplanes $H, \ldots , g^D \cdot H$ pairwise cross, which is impossible in the cube complex $X$ with dimension at most $D$.

As a consequence, there exist $0 \leq i < j \leq D$ such that $g^i \cdot H \cap g^j \cdot H = \emptyset$. Let $k = j-i \leq D$, we have $H \cap g^k \cdot H = \emptyset$. In particular, we have either $H^+ \subset g^k \cdot H^+$ or $g^k \cdot H^+ \subset H^+$. Since $k$ divides $D!$, we have either $H^+ \subset g^{D!} \cdot H^+$ or $g^{D!} \cdot H^+ \subset H^+$.

As $g^{D!} \cdot x \not\in g^{D!} \cdot H^+$ and $g^{D!} \cdot x \in H^+$, we conclude that $g^{D!} \cdot H^+ \subset H^+$.
\ep

We now prove a very similar statement, but with the weaker assumption that $x \in \Min(g^{D!})^{(0)}$ instead of $x \in \Min(g)^{(0)}$.

\blem \label{lem:case M(gD!)} Let $X$ be a CAT(0) cube complex of dimension at most $D$, and let $g$ be a combinatorially hyperbolic isometry of $X$. Then for any $x \in \Min(g^{D!})^{(0)}$, and for any hyperplane $H$ separating $x$ and $g^{D!} \cdot x$, $g^{D!}$ skewers $H$. \elem

\bp 
Fix $x \in \Min(g^{D!})^{(0)}$ and $y \in \Min(g)^{(0)}$. By contradiction, assume that there exists a hyperplane $H$ such that $x \in H^-$, $g^{D!} \cdot x \in H^+$ and $H$ is not skewered by $g^{D!}$. According to Lemma~\ref{lem:case M(g)}, for every $n,m \in \Z$, $H$ does not separate $g^{nD!} \cdot y$ and $g^{mD!} \cdot y$. By symmetry, assume that $\forall n \in \Z$, $g^{nD!} \cdot y \in H^+$.

As a consequence, for every $n \geq 0$, $g^{nD!} \cdot H$ separates $x$ and $y$. Since only $d_1(x,y)$ hyperplanes separate $x$ and $y$, we deduce that there exists $n>0$ such that $g^{nD!} \cdot H=H$. This contradicts the fact that a combinatorial geodesic from $x$ to $g^{nD!} \cdot x$ goes via $g^{D!} \cdot x$ as $x \in \Min(g^{D!})^{(0)}$ and crosses $H$, whereas $g^{nD!} \cdot H=H$ does not separate $x$ and $g^{nD!} \cdot x$.
\ep

\blem \label{lem:number of hyperplanes} Let $X$ be a CAT(0) cube complex of dimension at most $D$, and let $g$ be a combinatorially hyperbolic isometry of $X$ of translation length $\delta$. Then the set
$$\{\mbox{hyperplanes of $X$ skewered by $g^{D!}$}\}/<g^{D!}>$$ has cardinality $D!\delta$. \elem

\bp
Let $h=g^{D!}$, and fix $x \in \Min(h)^{(0)}$.

Let $H$ be a hyperplane skewered by $h$, so that $h \cdot H^+ \subset H^+$. Choose $n \in \Z$ such that $d_1(x,h^n \dot H)$ is minimal. Without loss of generality, assume that $x \in h^n \cdot H^+$. Since $h^{n+1} \cdot H$ is disjoint from $h^n \cdot H$ and $h^{n+1} \cdot H^+ \subset h^n \cdot H^+$, we deduce that $h^{n+1} \cdot H$ does not separate $h^n \cdot H$ and $x$, and so $h^{n+1} \cdot H$ separates $x$ and $h \cdot x$.

The number of hyperplanes separating $x$ and $h \cdot x$ is equal to $D!\delta$, so the cardinality of $\{$hyperplanes of $X$ skewered by $h\}/<h>$ is at most $D!\delta$.

\mk

Fix a hyperplane $H$ separating $x$ and $h \cdot x$: according to Lemma~\ref{lem:case M(gD!)}, $h$ skewers $H$. For instance, $x \in H^-$, $h \cdot x \in H^+$ and $h \cdot H^+ \subset H^+$. As a consequence, if $n>0$ then $x,h \cdot x \in h^n \cdot H^-$ so $h^n \cdot H$ does not separate $x$ and $h \cdot x$. Similarly, if $n<0$ then $x,h \cdot x \in h^n \cdot H^+$ so $h^n \cdot H$ does not separate $x$ and $h \cdot x$.

We conclude that the $D!\delta$ hyperplanes separating $x$ and $h \cdot x$ are disjoint $<h>$-orbits, hence the cardinality of $\{$hyperplanes of $X$ skewered by $h\}/<h>$ is exactly $D!\delta$.
\ep

\bpro \label{pro:M(g) median} Let $X$ be a CAT(0) cube complex of dimension at most $D$, and let $g$ be a combinatorially hyperbolic isometry of $X$. Then $\Min(g^{D!})^{(0)}$ is a median subalgebra of $X^{(0)}$, i.e. it is stable under the median of $X^{(0)}$. \epro

\bp 
Let $h=g^{D!}$, and let $\delta$ denote the translation length of $g$. Let $\mu$ denote the median of $X$.

\mk

Let $x,y,z \in \Min(h)^{(0)}$, and let $m=\mu(x,y,z) \in X^{(0)}$. Let $H$ be a hyperplane separating $m$ and $h \cdot m$, for instance $m \in H^-$ and $h \cdot m \in H^+$. Since $m=\mu(x,y,z) \in H^-$ which is convex, at least two vertices among $x$, $y$ and $z$ belong to $H^-$: we can assume that $x,y \in H^-$. Similarly, $h \cdot m=\mu(h \cdot x,h \cdot y,h \cdot z) \in H^+$ which is convex, at least two vertices among $h \cdot x$, $h \cdot y$ and $h \cdot z$ belong to $H^+$: we can assume that $h \cdot x,h \cdot z \in H^+$. As a consequence, $H$ separates $x$ and $h \cdot x$, so by Lemma~\ref{lem:case M(gD!)} $H$ is skewered by $h$.

According to Lemma~\ref{lem:number of hyperplanes}, we conclude that at most $D!\delta$ hyperplanes separate $m$ and $h \cdot m$. Since the translation length of $h$ is $D!\delta$, we conclude that $d(m,h \cdot m)=D!\delta$, so $m \in \Min(h)^{(0)}$. \ep

\brk There exists a combinatorially hyperbolic isometry $g$ of a locally finite CAT(0) cube complex such that for any $n \geq 1$, $\Min(g^n)$ is not convex.

For instance, consider a infinite, rooted, binary tree, where each edge is replaced by the diagonal of a square: this defines a CAT(0) square complex $T$. Then consider the isometry $g$ of $T$ defined recursively as shown in Figure~\ref{fig:recursive}. Notice that the fixed point set of $g^2$ in $T$ is equal to the union of the two diagonals of the squares adjacent to the root, and it is not convex for the combinatorial distance $d_1$. Furthermore, for any $n \in \N$, $g^{2^n}$ fixes some vertex $v$ (in fact, any vertex at the level $n$) of $T$ and acts on the subcomplex $T_v$ (corresponding to the subtree defined by $v$) as $g$. We deduce that, for any $n \geq 1$, the fixed point set of $g^{2^n}$ is not convex (we could argue similarly that, for every even $n \geq 1$, the fixed point set of $g^{2n}$ is not convex). 

If we want to find an example where $g$ is a combinatorially hyperbolic isometry, it is sufficient to consider the direct product $T \times \R$, where $g$ acts on $\R$ by a unit translation.
\erk

\begin{figure}
\begin{center}
\begin{tikzpicture}
\def \p {0.05}

\draw[fill] (-3,0) circle (0.1) node(a) {};
\draw[fill] (-2,0) circle (\p) node(b) {};
\draw[fill] (-1.5,1) circle (\p) node(c) {};
\draw[fill] (-2.5,1) circle (\p) node(d) {};
\draw[fill] (-4,0) circle (\p) node(-b) {};
\draw[fill] (-4.5,1) circle (\p) node(-c) {};
\draw[fill] (-3.5,1) circle (\p) node(-d) {};
\draw[fill] (-1.5,1)+(80:1) circle (\p) node(c1) {};
\draw[fill] (-1.5,1)+(60:1) circle (\p) node(c2) {};
\draw[fill] (-1.5,1)+(70:1.8) circle (\p) node(c3) {};
\draw[fill] (-1.5,1)+(100:1) circle (\p) node(c4) {};
\draw[fill] (-1.5,1)+(120:1) circle (\p) node(c5) {};
\draw[fill] (-1.5,1)+(110:1.8) circle (\p) node(c6) {};
\draw[fill] (-4.5,1)+(80:1) circle (\p) node(-c1) {};
\draw[fill] (-4.5,1)+(60:1) circle (\p) node(-c2) {};
\draw[fill] (-4.5,1)+(70:1.8) circle (\p) node(-c3) {};
\draw[fill] (-4.5,1)+(100:1) circle (\p) node(-c4) {};
\draw[fill] (-4.5,1)+(120:1) circle (\p) node(-c5) {};
\draw[fill] (-4.5,1)+(110:1.8) circle (\p) node(-c6) {};

\draw (a.center) -- (b.center) -- (c.center) -- (d.center) -- (a.center);
\draw (a.center) -- (-b.center) -- (-c.center) -- (-d.center) -- (a.center);
\draw [->>,>=triangle 60, color=blue] (b.center) -- (d.center);
\draw [->,>=triangle 60, color=red] (-b.center) -- (-d.center);
\draw (c.center) -- (c1.center) -- (c3.center) -- (c2.center) -- (c.center);
\draw (c.center) -- (c4.center) -- (c6.center) -- (c5.center) -- (c.center);
\draw (-c.center) -- (-c1.center) -- (-c3.center) -- (-c2.center) -- (-c.center);
\draw (-c.center) -- (-c4.center) -- (-c6.center) -- (-c5.center) -- (-c.center);

\draw[fill] (3,0) circle (0.1) node(a') {};
\draw[fill] (2,0) circle (\p) node(b') {};
\draw[fill] (1.5,1) circle (\p) node(c') {};
\draw[fill] (2.5,1) circle (\p) node(d') {};
\draw[fill] (4,0) circle (\p) node(-b') {};
\draw[fill] (4.5,1) circle (\p) node(-c') {};
\draw[fill] (3.5,1) circle (\p) node(-d') {};
\draw[fill] (1.5,1)+(80:1) circle (\p) node(c1') {};
\draw[fill] (1.5,1)+(60:1) circle (\p) node(c2') {};
\draw[fill] (1.5,1)+(70:1.8) circle (\p) node(c3') {};
\draw[fill] (1.5,1)+(100:1) circle (\p) node(c4') {};
\draw[fill] (1.5,1)+(120:1) circle (\p) node(c5') {};
\draw[fill] (1.5,1)+(110:1.8) circle (\p) node(c6') {};
\draw[fill] (4.5,1)+(80:1) circle (\p) node(-c1') {};
\draw[fill] (4.5,1)+(60:1) circle (\p) node(-c2') {};
\draw[fill] (4.5,1)+(70:1.8) circle (\p) node(-c3') {};
\draw[fill] (4.5,1)+(100:1) circle (\p) node(-c4') {};
\draw[fill] (4.5,1)+(120:1) circle (\p) node(-c5') {};
\draw[fill] (4.5,1)+(110:1.8) circle (\p) node(-c6') {};
\node (gg) at (-4.5,3) {};
\node (gd) at (-1.5,3) {};
\node (dg) at (1.5,3) {};
\node (dd) at (4.5,3) {};

\draw (a'.center) -- (b'.center) -- (c'.center) -- (d'.center) -- (a'.center);
\draw (a'.center) -- (-b'.center) -- (-c'.center) -- (-d'.center) -- (a'.center);
\draw [->>,>=triangle 60, color=blue] (d'.center) -- (b'.center);
\draw [->,>=triangle 60, color=red] (-b'.center) -- (-d'.center);
\draw (c'.center) -- (c1'.center) -- (c3'.center) -- (c2'.center) -- (c'.center);
\draw (c'.center) -- (c4'.center) -- (c6'.center) -- (c5'.center) -- (c'.center);
\draw (-c'.center) -- (-c1'.center) -- (-c3'.center) -- (-c2'.center) -- (-c'.center);
\draw (-c'.center) -- (-c4'.center) -- (-c6'.center) -- (-c5'.center) -- (-c'.center);

\draw [->,>=triangle 60,auto] (a) edge [bend right] node {g} (a');
\draw [->,>=triangle 60,auto] (gg) edge [bend left] node {id} (dd);
\draw [->,>=triangle 60,auto] (gd) edge [bend left] node {g} (dg);

\end{tikzpicture}
\end{center}
\caption{A pathological isometry of a CAT(0) square complex.}
\label{fig:recursive}
\end{figure}

We will now show that the minimal set of an abelian group is not empty.

\bpro \label{pro:minset_a}
Let $X$ be a CAT(0) cube complex of dimension at most $D$, and let $A$ be a finitely generated abelian group of cubical automorphisms of $X$. Then $\Min(A^{2D!})^{(0)} = \bigcap_{a \in A} \Min(a^{2D!})^{(0)}$ is not empty.
\epro

\bp
We will prove the result by induction on the numbers of generators of $A$. When $A$ is cyclic it follows from Proposition~\ref{pro:min_vertex}. Assume that $A=\<B,g\>$, where $\Min(B^{2D!})^{(0)} \neq \emptyset$.

Since $g$ commutes with $B^{2D!}$, $g$ preserves $\Min(B^{2D!})^{(0)}$. According to Proposition~\ref{pro:M(g) median}, $\Min(B^{2D!})^{(0)}$ is a median subalgebra of $X^{(0)}$. The element $g$ acts as a cubical automorphism on the CAT(0) cube complex $Y=X(\Min(B^{2D!})^{(0)})$. Since $\Min(B^{2D!})^{(0)}$ is a median subalgebra of $X^{(0)}$ which has rank at most $D$, we deduce that $Y$ has dimension at most $D$.

According to Proposition~\ref{pro:min_vertex}, there exists a vertex $x \in Y^{(0)} = \Min(B^{2D!})^{(0)}$ which lies in $\Min_Y(g^{2D!})$ (the minimal set of $g^{2D!}$ considered as an isometry of $Y$). This means that $\{g^{2nD!} \cdot x,n \in \Z\}$ lie on a combinatorial geodesic in $Y$. Since $\Min(B^{2D!})^{(0)}$ is a median subalgebra of $X^{(0)}$, we deduce that $\{g^{2nD!} \cdot x,n \in \Z\}$ lie on a combinatorial geodesic in $X$. This precisely means that $x \in \Min(g^{2D!})$. This vertex $x$ of $X$ belongs to $\Min(g^{2D!}) \cap \Min(B^{2D!}) = \Min(A^{2D!})$, which proves that $\Min(A^{2D!})^{(0)} \neq \emptyset$.
\ep

\bpro \label{pro:min set cocompact} Let $G$ be a group acting geometrically on a locally finite CAT(0) cube complex $X$ of dimension at most $D$, and let $A$ be an abelian subgroup of $G$ consisting of $2D!^\text{th}$ powers of combinatorially hyperbolic isometries. Then the centralizer $Z_G(A)$ of $A$ in $G$ has finite index in the normalizer $N_G(A)$ of $A$ in $G$, and $N_G(A)$ acts geometrically on $\Min(A)^{(0)}$. \epro

\bp According to Proposition~\ref{pro:minset_a}, the minimal set $\Min(A)^{(0)}$ is not empty.

\mk

The action of $N_G(A)$ on $X$ is proper and stabilizes $\Min(A)^{(0)}$, so it induces a proper action on $\Min(A)^{(0)}$.

Assume that the action of $Z_G(A)$ on $\Min(A)^{(0)}$ is not cocompact: since $G$ acts properly and cocompactly on $X$, there exist $C \geq 0$, $x \in \Min(A)^{(0)}$ and $(h_n)_{n \in \N} \in G^\N$ such that $\forall n \in \N, d(h_n \cdot x,\Min(A)^{(0)}) \leq C$ and the cosets $(h_nZ_G(A))_{n \in \N} \in \left(G/Z_G(A)\right)^\N$ are pairwise distinct.

\mk

According to the flat torus theorem (see~\cite[Theorem~7.1]{bridson_haefliger}), the abelian group $A$ acts properly by semisimple isometries on the CAT(0) space $X$, so $A$ is finitely generated. Fix $a_1,\dots,a_r$ some generators of $A$. Fix some $1 \leq i \leq r$.

Let $\delta_i$ denote the combinatorial translation length of $a_i$. For all $n \in \N$, since $d(h_n\cdot x,\Min(a_i)^{(0)}) \leq C$ we have $d(h_n \cdot x, a_ih_n \cdot x) \leq \delta_i+2C$. So $d(x,h_n^{-1}a_ih_n \cdot x) \leq \delta_i+2C$. Since $X$ is locally finite, up to passing to a subsequence, we may assume that $\forall n,m \in \N, h_n^{-1}a_ih_n \cdot x = h_m^{-1}a_ih_m \cdot x$. Since the action of $G$ on $X$ is proper, the stabilizer of $x$ is finite, so up to passing to a new subsequence, we may assume that $\forall n,m \in \N, h_n^{-1}a_ih_n = h_m^{-1}a_ih_m$. So $\forall n,m \in \N, h_nh_m^{-1}$ centralizes $a_i$.

\mk

If we apply this for every $1 \leq i \leq r$, we obtain up to passing to a new subsequence that $\forall n,m \in \N, h_nh_m^{-1}$ centralizes $a_1,\dots,a_r$. Since $a_1,\dots,a_r$ span $A$, we deduce that $\forall n,m \in \N, h_nh_m^{-1} \in Z_G(A)$. This contradicts the assumption that the cosets $(h_nZ_G(A))_{n \in \N} \in \left(G/Z_G(A)\right)^\N$ are pairwise distinct.

As a consequence, the induced action of $Z_G(A)$ on $\Min(A)^{(0)}$ is proper and cocompact. Since the action of $N_G(A)$ on $\Min(A)^{(0)}$ is also proper, we deduce that $Z_G(A)$ has finite index in $N_G(A)$.
 \ep

We obtain now the proof of Theorem~\ref{thm:cubulation of multiple centralizer}.

\bp Let $a_1,\dots,a_r$ be some generators of $A$. For each $1 \leq i \leq r$, by Proposition~\ref{pro:M(g) median} $\Min(a_i)^{(0)}$ is a median subalgebra of $X^{(0)}$. As a consequence, $\Min(A)^{(0)}=\bigcap_{i=1}^r \Min(a_i)^{(0)}$ is a also a median subalgebra of $X^{(0)}$. According to Proposition~\ref{pro:minset_a}, it is not empty. By Proposition~\ref{pro:min set cocompact}, $N_G(A)$ acts properly cocompactly on $\Min(A)^{(0)}$. Theorem~\ref{thm:CCC from median} concludes the proof.

Note that the CAT(0) cube complex $X(\Min(A)^{(0)})$ has dimension at most $D$, since $\Min(A)^{(0)}$ is a median subalgebra of $X^{(0)}$ which has rank at most $D$.
\ep

\brk Note that the distances induced on $\Min(A)^{(0)}$ by $X$ and by $X(\Min(A)^{(0)})$ may be different.

For instance, consider the action of $A=\<g\> \simeq \Z$ on $X=\R^2$ by $g \cdot (x,y) = (y+1,x+1)$. Then $\Min(g)^{(0)}=\{(n,n) \st n \in \Z\}$. Inside $X$, the combinatorial distance between $(n,n)$ and $(n+1,n+1)$ is equal to $2$. But as a median algebra, $\Min(g)^{(0)}$ is isomorphic to $\Z$, and hence $X(\Min(g)^{(0)})$ is isomorphic to the standard cubical structure on $\R$. So in $X(\Min(g)^{(0)})$, the combinatorial distance between $(n,n)$ and $(n+1,n+1)$ is equal to $1$. This is because the two hyperplanes of $\R^2$ separating $(n,n)$ and $(n+1,n+1)$ define the same partition of $\Min(g)^{(0)}$.

This example does not satisfy the assumption of Theorem~\ref{thm:cubulation of multiple centralizer} that elements of $A$ are fourth powers of isometries, but nevertheless it illustrates the difference between the two cubical structures on $\Min(A)^{(0)}$.
\erk

\section{Convex-cocompact subgroups}

In this article, the \emph{rank} of a finitely generated abelian or virtually abelian group is the minimal number of generators of a finite index free abelian subgroup.

\bdf A subgroup $A$ of a group $G$ acting geometrically on a CAT(0) cube complex $X$ is said to be \emph{convex-cocompact in $X$} if there exists a convex subcomplex $Y \subset X$ which is $A$-invariant and $A$-cocompact. \edf

We now give an equivalent characterization of convex-cocompact subgroups, starting with a small Lemma.

\blem \label{lem:1neighbourhood}
Let $X$ be a finite-dimensional cube complex and let $Y \subset X$ be a convex subcomplex. For any $n \in \N$, there exists a convex subcomplex $Y_n$ of $X$ containing the combinatorial neighbourhood of radius $n$ of $Y$, and such that $Y_n$ is contained in a bounded neighbourhood of $Y$. Furthermore, if $g$ is an automorphism of $X$ preserving $Y$, then $g$ preserves $Y_n$ for every $n \in \N$.
\elem

\bp
It is sufficient to prove the statement for $n=1$. Let $Z$ denote the full subcomplex of $X$ whose vertices are all vertices $x \in X^{(0)}$ such that any two hyperplanes separating $x$ and $Y$ cross. Then $Z$ contains the combinatorial $1$-neighbourhood of $Y$, we will prove that $Z$ is convex.

Fix $z,z' \in Z^{(0)}$, and consider a vertex $x$ on a combinatorial geodesic from $z$ to $z'$, we will prove that $x \in Z$. Assume that there exist two distinct hyperplanes $H,H'$ separating $x$ and $Y$. We will prove that $H$ and $H'$ cross. Since $x$ lies on a combinatorial geodesic from $z$ to $z'$, we deduce that $H$ separates $z$ and $Y$, or separates $z'$ and $Y$. We deduce the same for $H'$.

If both $H$ and $H'$ separate $z$ and $Y$, then $H$ and $H'$ cross since $z \in Z$.

So we can assume that $H$ separates $\{z,x\}$ and $\{z'\} \cup Y$, and that $H'$ separates $\{z',x\}$ and $\{z\} \cup Y$. Hence the four intersections of the half-spaces of $H$ and $H'$ are non-empty, so $H$ and $H'$ cross.

As a consequence, $Z$ is convex.

Finally, for any vertex $z \in Z^{(0)}$, let $y \in Y^{(0)}$ denote the gate projection of $z$ onto $Y$. Then by definition of $Z$, all hyperplanes separating $z$ and $y$ pairwise cross. As a consequence, $z$ and $y$ belong to a common cube of $X$. We deduce that $y$ and $z$ are at distance at most the dimension of $X$. In particular, $Z$ is contained in a bounded neighbourhood of $Y$.\ep

\bpro \label{pro:cvx_cocompact_some_every} Let $A$ be a subgroup of a group $G$ acting geometrically on a CAT(0) cube complex $X$. Then $A$ is convex-cocompact if and only if for every $x \in X^{(0)}$ (equivalently, for some $x \in X^{(0)}$), $A$ acts cocompactly on $\Hull(A \cdot x)$.
\epro

\bp
If $A$ is convex-cocompact, for any $x \in Y$, we have $\Hull(A \cdot x) \subset Y$ and so $A$ acts cocompactly on $\Hull(A \cdot x)$.

Conversely, assume that there exists a vertex $x \in X^{(0)}$ such that $A$ acts cocompactly on $Y=\Hull(A \cdot x)$. For any vertex $y \in X^{(0)}$, according to Lemma~\ref{lem:1neighbourhood}, there exists a convex subcomplex $Z$ of $X$ containing $Z$ and $y$, contained in a neighbourhood of $Y$, which is furthermore $A$-invariant. Since $X$ is locally finite, we deduce that $A$ acts cocompactly on $Z$. Since $\Hull(A \cdot y) \subset Z$, we conclude that $A$ acts cocompactly on $\Hull(A \cdot y)$. 
\ep

\brk Note that being convex-cocompact depends on the CAT(0) cube complex $X$: see for instance Subsection~\ref{subsec:cubulation_even}. \erk

We will now state a few technical results about convex-cocompact subgroups which will be used in the sequel.

\blem \label{lem:intersection of convex subgroups} Let $G$ be a group acting geometrically on a CAT(0) cube complex $X$, and let $A,B$ be subgroups of $G$ which are convex-cocompact in $X$. Then $A \cap B$ is convex-cocompact in $X$. \elem

\bp Fix a vertex $x \in X^{(0)}$, and consider a sequence $(x_n)_{n \in \N}$ in $\Hull(A \cap B \cdot x)$. According to Proposition~\ref{pro:cvx_cocompact_some_every}, $A$ and $B$ act cocompactly on $\Hull(A \cdot x)$ and $\Hull(B \cdot x)$ respectively, so there exist $C>0$ and sequences $(a_n)_{n \in \N}$ in $A$ and $(b_n)_{n \in \N}$ in $B$ such that $\forall n \in \N, d(a_n \cdot x, x_n) \leq C$ and $d(b_n \cdot x,x_n) \leq C$. As a consequence, $\forall n \in \N, d(b_n^{-1}a_n \cdot x,x) \leq 2C$. Since $G$ acts properly on the CAT(0) cube complex, we deduce that, up to passing to subsequences, we have $\forall n,m \in \N, b_m^{-1}a_m=b_n^{-1}a_n$, so $a_na_m^{-1} = b_nb_m^{-1} \in A \cap B$. So for all $n \in \N$, we have $d(a_na_0^{-1} \cdot x,x_n) \leq d(a_n \cdot x,x_n)+d(a_na_0^{-1} \cdot x,a_n \cdot x) \leq C+d(a_0^{-1} \cdot x,x)$ is bounded. Since $\forall n \in \N, a_na_0^{-1} \in A \cap B$, this proves that $A \cap B$ acts cocompactly on $\Hull(A \cap B \cdot x)$. Hence $A \cap B$ is convex-cocompact in $X$.
\ep

\bdf A virtually abelian subgroup $A$ of a group $G$ is called \emph{highest} if for any virtually abelian subgroup $B$ of $G$ such that $A \cap B$ has finite index in $A$, we have that $A \cap B$ has finite index in $B$. \edf

We now recall the following recent result from D.~Wise and D.~Woodhouse.

\bthm[Cubical flat torus theorem \cite{wise_woodhouse}] \label{thm:flat_torus} Let $G$ be a group acting geometrically on a CAT(0) cube complex $X$. Let $A$ be a highest virtually abelian subgroup of $G$. Then $A$ is convex-cocompact in $X$. \ethm

We can now prove the following, which is the main technical part in the proof of Theorem~\ref{thm:stability_central_quotients} stated in the introduction.

\blem \label{lem:splitting} Let $G$ be a group acting geometrically, combinatorially semisimply by isometries on a CAT(0) cube complex $X$, median minimally. Let $W$ be a central subgroup of $G$ which is convex-cocompact in $X$. Then $X$ splits as a product of two convex cube subcomplexes $X \simeq Y \times Z$, where $G$ preserves this splitting, $W$ acts with finite index kernel $W'$ on $Y$ and $G/W'$ acts geometrically on $Y$, with $\dim Y \leq \dim X - \rk W$. Furthermore, if $A \subset G$ is convex-cocompact in $X$, then $AW'/W'$ is convex-cocompact in $Y$. \elem

\bp Let $D$ denote the dimension of $X$, and let $W'=W^{2D!^2}=\< w^{2D!^2} \st w \in W\>$. Since $W$ is abelian and finitely generated, $W'$ has finite index in $W$. 

Let $Skew(W')$ denote the set of hyperplanes of $X$ that are skewered by at least one element of $W'$. Let $Stab(W')$ denote the set of hyperplanes of $X$ that are stabilized by all elements of $W'$. We will prove that the set of hyperplanes of $X$ is the disjoint $G$-invariant union $Skew(W') \sqcup Stab(W')$. By definition, the sets $Skew(W')$ and $Stab(W')$ are disjoint. Since $W$ is central in $G$, both sets $Skew(W')$ and $Stab(W')$ are $G$-invariant.

\mk

Let us first remark that each hyperplane in $Skew(W')$ intersects each hyperplane in $Stab(W')$. Fix $H \in Skew(W')$, it is skewered by some $w \in W'$. Fix $H' \in Stab(W')$, it is stabilized by $w$. Since $H$ is skewered by $w$, any $\<w\>$-orbit in $X$ intersects both half-spaces of $H$. In particular, if $x \in H'$, we know that $\<w\> \cdot x$ lies in $H'$ and intersects both half-spaces of $H$. Hence $H$ intersects $H'$.

\mk

So we have to prove that every hyperplane of $X$ is in the union $Skew(W') \sqcup Stab(W')$. Fix a hyperplane $H$ which is not in $Skew(W')$, we will prove that $H \in Stab(W')$.

\mk

Since $W$ is abelian, according to Proposition~\ref{pro:minset_a}, we have $\Min(W^{2D!})^{(0)} \neq \emptyset$. According to Proposition~\ref{pro:M(g) median}, $\Min(W^{2D!})^{(0)}$ is a median subalgebra of $X^{(0)}$. As $W^{2D!}$ is normal in $G$, we deduce that $G$ stabilizes $\Min(W^{2D!})^{(0)}$. As the action of $G$ on $X^{(0)}$ is median minimal, we deduce that $X^{(0)}=\Min(W^{2D!})^{(0)}$.

\mk

Note that any decreasing sequence of $W^{2D!}$-cocompact subcomplexes of $X$ terminates, so there exists a vertex $x \in X^{(0)}$ for which the action of $W^{2D!}$ on $Z=\Hull(W^{2D!} \cdot x)$ is such that $Z$ is the smallest non-empty convex $W^{2D!}$-invariant subcomplex of $Z$.

\mk

Any hyperplane intersecting $Z$ separates two points in $W^{2D!} \cdot x$, so according to Lemma~\ref{lem:case M(g)} such a hyperplane is skewered by an element of $W^{2D!}$. Hence $H$ does not intersect $Z$.

\mk

Since the action of $G$ on $X$ is median minimal, the orbit $G \cdot x$ is not contained in a single half-space of $H$, so there exists $g \in G$ such that $H$ separates $x$ and $g \cdot x$. Say $Z \subset H^-$ and $g \cdot x \in H^+$.

\mk

We will show that for every $w \in W^{2D!}$, we have $w^{D!} \cdot H=H$. Fix $w \in W^{2D!} \bs \{1\}$. We will first show that $H$ separates $Z$ and $g\<w\> \cdot x$. So we want to show that $g \<w\> \cdot x \subset H^+$.

\medskip

By contradiction, assume that there exists $n_0 \neq 0$ such that $gw^{n_0} \cdot x \in H^-$. Without loss of generality (up to replacing $w$ with $w^{-1}$), we can assume that $n_0>0$.

As $g \cdot x \in X^{(0)} = \Min(W^{2D!})^{(0)}$, we deduce that $(gw^k \cdot x)_{k \in \Z}$ lies on a combinatorial geodesic of $X$. As a consequence, for every $n \leq 0$, we have $gw^n \cdot x \in H^+$. So for every $n \leq 0$, we have $g \cdot x \in w^{-n} \cdot H^+$ and $x \in Z = w^{-n} \cdot Z \subset w^{-n} \cdot H^-$.

Since finitely many hyperplanes separate $x$ and $g \cdot x$, we conclude that there exists $n<0$ such that $w^n \cdot H = H$. On one hand, as $|n|n_0 \geq n_0$ and $(gw^k \cdot x)_{k \in \Z}$ lies on a combinatorial geodesic of $X$, we have $gw^{|n|n_0} \cdot x \in H^-$. On the other hand, as $w^{|n|n_0} \cdot H^+ = H^+$ and $g \cdot x \in H^+$ we have $gw^{|n|n_0} \cdot x \in H^+$. This is a contradiction.

\medskip

As a consequence, every hyperplane separating $Z$ and $g \cdot x$ separates $Z$ and $g\<w\> \cdot x$. Let ${\cal H}$ denote the finite set of hyperplanes separating $Z$ and $g \cdot x$. We have seen that $w$ preserves the set ${\cal H}$, and acts as a bijection $\sigma$ on ${\cal H}$. Since $w$ preserves $Z$, we know that hyperplanes of ${\cal H}$ in the same $\sigma$-orbit are at the same distance of $Z$, and thus cannot be nested. Hence each $k$-cycle of $\sigma$ in ${\cal H}$ corresponds to $k$ pairwise crossing hyperplanes, so $k \leq D$. As a consequence, $\sigma^{D!}=1$. Since $H \in {\cal H}$, we conclude that $w^{D!} \cdot H=H$. Furthermore, since every element of $G$ acts combinatorially semisimply, we deduce that $w^{D!}$ stabilizes each half-space of $H$.

\mk

We have proved that for every $w \in W^{2D!}$, we have $w^{D!} \cdot H=H$, so $H \in Stab(W')$.

\mk

The $G$-equivariant disjoint decomposition $Skew(W') \sqcup Stab(W')$ defines a $G$-equivariant isomorphism $X \simeq Y \times Z$, where $Y$ is the cube complex dual to the set of hyperplanes $Stab(W')$. Note that $Z$ is also isomorphic to the cube complex dual to the set of hyperplanes $Skew(W')$.

\mk

Since $W$ acts properly on $Z$, we have $\dim Z \geq \rk W$, so $\dim Y \leq \dim X - \rk W$. Furthermore, by definition of $Stab(W')$, we know that $W'$ acts trivially on $Y$, so that the $G$-action on $Y$ factors through $G/W'$.

As $G$ acts cocompactly on $X$, we deduce that $G/W'$ acts cocompactly on $Y$. Since $W'$ acts properly on $Z$, we deduce that $G/W'$ acts properly on $Y$. In conclusion, $G/W'$ acts geometrically on $Y$.

\mk

Furthermore, if $A \subset G$ is convex-cocompact in $X$, then $A$ acts cocompactly on a convex subcomplex $M_A$ of $X$. Since $X \simeq Y \times Z$, the convex subcomplex $M_A$ splits as $M_A \simeq M_{A,Y} \times M_{A,Z}$, such that $AW'/W'$ acts geometrically on the convex subcomplex $M_{A,Y}$ of $Y$.
\ep

\blem \label{lem:convex_cocompact_intersection_centralizer}
Let $G$ be a group acting geometrically on a CAT(0) cube complex $X$. Assume that $H \subset G$ acts geometrically on a median subalgebra $M \subset X^{(0)}$, with associated CAT(0) cube complex $X(M)$. Assume that $A$ is a convex-cocompact subgroup of $G$. Then $A \cap H$ is convex-cocompact in $X(M)$.
\elem

\bp
Consider a convex subcomplex $Y$ of $X$ such that $A$ acts properly and cocompactly on $Y$. Up to considering some convex neighbourhood of $Y$, we may assume that $Y$ and $M$ intersect. Let $x_0 \in Y \cap M$.

The group $A \cap H$ acts properly on $Y \cap M$. We will prove that the action is cocompact. By contradiction, assume that there exists a sequence $(y_n)_{n \in \N}$ in $Y \cap M$ such that $d_1(y_n, A \cap H \cdot x_0) \ral{n \ra \pif} \pif$. Since the action of $G$ on $X$ is cocompact, there exists $D \geq 0$ such that for each $n \in \N$, there exists $g_n \in G$ such that $d_1(y_n,g_n \cdot x_0) \leq D$, so that $d_1(g_n \cdot x_0, A \cap H \cdot x_0) \ral{n \ra \pif} \pif$.

\mk

Since the action of $A$ on $Y$ is cocompact and $d_1(g_n \cdot x_0,Y) \leq D$ for all $n \in \N$, there exists a finite subset $K$ of $G$ such that, for each $n \in \N$, we have $g_n \in AK$. Similarly, since the action of $H$ on $M$ is cocompact and $d_1(g_n \cdot x_0,M) \leq D$ for all $n \in \N$, there exists a finite subset $K'$ of $G$ such that, for each $n \in \N$, we have $g_n \in HK'$. Up to passing to a subsequence, we may assume that there exists $k \in K$ and $k' \in K'$ such that for each $n \in \N$, we have $g_n \in Ak \cap Hk'$. In particular, for each $n \in \N$ we have $g_n g_0^{-1} \in A \cap H$. So $d_1(g_n \cdot x_0, A \cap H \cdot x_0) = d_1((g_n g_0^{-1})g_0 \cdot x_0, A \cap H \cdot x_0)=d_1(g_0 \cdot x_0, A \cap H \cdot x_0)$ is bounded, which is a contradiction.

\mk

So we have proved that the group $A \cap H$ acts properly and cocompactly on $Y \cap M$. Since $Y^{(0)}$ is convex in the median algebra $X^{(0)}$, according to Lemma~\ref{lem:median inter convex}, we deduce that $M \cap Y^{(0)}$ is convex in the median subalgebra $M$. This implies that $A \cap H$ is a convex-cocompact subgroup in $X(M)$.
\ep

\section{Non-cubicality criterion}

We will now summarize the main two stability results for virtual cubulation that we will use. They are slightly more precise than Theorem~\ref{thm:stability_normalizers} (only the part about centralizers) and Theorem~\ref{thm:stability_central_quotients}.

\bpro \label{pro:cubulation_centralizer}
Let $G$ be a group acting geometrically on a CAT(0) cube complex $X$, and let $C$ be an abelian subgroup of $G$. There exists a finite index subgroup $C_0$ of $C$ such that the centralizer $Z_G(C_0)$ acts geometrically on a CAT(0) cube complex $Y$ with $\dim Y \leq \dim X$.

Furthermore, if $A \subset G$ is convex-cocompact in $X$, then $A \cap Z_G(C_0)$ is convex-cocompact in $Y$.
\epro

\bp
This is a consequence of Theorem~\ref{thm:cubulation of multiple centralizer} and Lemma~\ref{lem:convex_cocompact_intersection_centralizer}.
\ep

\bpro \label{pro:cubulation_central_quotient}
Let $G$ be a group acting geometrically on a CAT(0) cube complex $X$, and let $W$ be a central subgroup of $G$ which is convex-cocompact in $X$. Then $W$ has a finite index subgroup $W'$ such that $G/W'$ acts geometrically on a CAT(0) cube complex $Y$ with $\dim Y \leq \dim X - \rk W$.

Furthermore, if $A \subset G$ is convex-cocompact in $X$, then $AW'/W'$ is convex-cocompact in $Y$.
\epro

\bp
This is contained in Lemma~\ref{lem:splitting}. The combinatorially semisimple assumption follows from Theorem~\ref{thm:haglund}, up to passing to the cubical subdivision of $X$. Furthermore, one can always assume that the action is median minimal, up to passing to a smaller median subalgebra. Note that this does not increase the dimension of $X$: indeed if $M'$ is a median subalgebra of $M$, then cubes in $M'$ are cubes in $M$.\ep

We now give a slightly more general version of a result from Crisp and Paoluzzi (see~\cite{crisp_paoluzzi}), which studies proper semisimple actions of $B_3$ and $B_4$ on CAT(0) spaces. Note that there is no cocompactness assumption in this result, nor a CAT(0) cube complex.

\bpro[Crisp-Paoluzzi] \label{pro:proper action of Im} Let $p \in \N_{\geq 3}$ and consider the dihedral Artin-Tits group $A=A(p)=\<a,b \mid w_p(a,b)=w_p(b,a)\>$. Assume $A$ acts properly, by semisimple isometries on a CAT(0) space $X$. Then $a$, $z_{ab}$ and $b$ act by hyperbolic isometries, whose attracting endpoints in the visual boundary $\bord X$ are denoted $a(\pif)$, $z_{ab}(\pif)$ and $b(\pif)$. Furthermore, if we denote by $\sphericalangle$ the visual distance on $\partial X$, we have:
\bit
\item If $p$ is odd, then we have $\sphericalangle(a(\pif),z_{ab}(\pif)) < \f{\pi}{2}$ and $\sphericalangle(b(\pif),z_{ab}(\pif)) < \f{\pi}{2}$.
\item If $p$ is even, then we have $\sphericalangle(a(\pif),z_{ab}(\pif)) < \f{\pi}{2}$ or $\sphericalangle(b(\pif),z_{ab}(\pif)) < \f{\pi}{2}$.
\eit\epro

\bp We adapt here the proof of \cite[Theorem~4]{crisp_paoluzzi}. Without loss of generality, we may assume that the action of $A$ on $X$ is minimal. By properness, every infinite order element of $A$ acts by a hyperbolic isometry, in particular $a$, $b$ and $z_{ab}$. Then by \cite[Theorem~II.6.8]{bridson_haefliger}, $X$ is isometric to the product $\R \times Y$, where $Y$ is a CAT(0) space, and $Z(A)=\<z_{ab}\>$ acts by translation on $\R$ and trivially on $Y$. Let $\delta \in \R \bs \{0\}$ such that $z_{ab}$ acts on the $\R$ factor by a translation of $\delta$. Up to the choice of the orientation of $\R$, we may assume that $\delta>0$.

Since $a$ and $b$ commute with $z_{ab}$, they preserve the decomposition $X \simeq \R \times Y$ and preserve the orientation of the $\R$ factor. In particular, let $\alpha,\beta \in \R$ such that $a$ and $b$ act on the $\R$ factor by translations of $\alpha$ and $\beta$ respectively.

\bit
\item If $p$ is odd, then $a$ and $b$ are conjugated by $w_p(a,b)$ in $A$, we deduce that $\alpha=\beta$. But $z_{ab}=w_{2p}(a,b)$, so we have $\delta=2p\alpha$. As a consequence, we have $\alpha=\beta>0$. This implies that the attracting endpoints of $a$ and $z_{ab}$ in $\bord X$ satisfy $\sphericalangle(a(\pif),z_{ab}(\pif)) < \f{\pi}{2}$.
\item If $p$ is even, then since $z_{ab}=w_{p}(a,b)$, we deduce that $p\alpha+p\beta=\delta>0$. As a consequence, $\alpha>0$ or $\beta>0$. This implies that $\sphericalangle(a(\pif),z_{ab}(\pif)) < \f{\pi}{2}$ or $\sphericalangle(b(\pif),z_{ab}(\pif)) < \f{\pi}{2}$.
\eit
\ep

\bpro \label{pro:2 abelian subgroups} Let $G$ be a group acting geometrically on a CAT(0) cube complex $X$, and let $A$, $B$ be subgroups of $G$ which are convex-cocompact in $X$, such that $A \cap B$ is finite. Then for each $a \in A$, $b \in B$ of infinite order, their attractive endpoints in $\bord X$ satisfy $\sphericalangle(a(\pif),b(\pif)) \geq \f{\pi}{2}$. \epro

\bp Let $M_A,M_B$ denote convex cube subcomplexes of $X$ on which $A,B$ respectively act geometrically.

Fix $x \in X^{(0)}$, and let $R \geq 0$ such that $d_1(x,M_A) \leq R$ and $d_1(x,M_B) \leq R$. Let $x_A \in M_A$ and $x_B \in M_B$ such that $d_1(x,x_A) \leq R$ and $d_1(x,x_B) \leq R$. Define
$$ S=\{y \in X^{(0)} \st d_1(y,M_A) \leq R \mbox{ and }d_1(y,M_B) \leq R\}.$$

We have $x \in S$. We claim that $S$ is finite: if not, since $X$ is locally compact, we can consider a sequence $(s_n)_{n \in \N}$ of pairwise distinct elements of $S$. Since $A$ and $B$ act geometrically on $M_A$ and $M_B$ respectively, up to passing to a subsequence, we deduce that there exist vertices $y_A \in M_A$, $y_B \in M_B$ and sequences $(a_n)_{n \in \N}$, $(b_n)_{n \in \N}$ of pairwise distinct elements in $A$ and $B$ respectively, such that the sequence $d_1(a_n \cdot y_A,b_n \cdot y_B)$ is bounded above. Since the action of $G$ on $X$ is proper and $X$ is locally compact, we can assume up to passing to a subsequence that the sequence $(b_n^{-1}a_n)_{n \in \N}$ is constant, hence for all $m,n \in \N$ we have $a_na_m^{-1}=b_nb_m^{-1} \in A \cap B$. As $A \cap B$ is finite, this is a contradiction. So $S$ is finite.

\medskip

From now on, fix $a \in A$ and $b \in B$ of infinite order. We will show that their attractive endpoints in $\bord X$ satisfy $\sphericalangle(a(\pif),b(\pif)) \geq \f{\pi}{2}$.

Let $\mu : X^3 \ra X$ denote the median on $X$. Fix $(\alpha_n)_{n \in \N},(\beta_n)_{n \in \N}$ sequences in $M_A^{(0)}$ (resp. $M_B^{(0)}$) converging to $a(\pif)$ (resp. $b(\pif)$). For each $n \in \N$, define $m_n = \mu(\alpha_n,\beta_n,x)$. As $\alpha_n$, $\beta_n$ and $x$ are vertices of $X$, $m_n$ is also a vertex of $X$.

Since $\mu$ is $1$-Lipschitz with respect to $d_1$ (see Lemma~\ref{lem:lipschitz}), we deduce that $d_1(m_n,M_A) \leq d_1(x,x_A)+d_1(\mu(\alpha_n,\beta_n,x_A),M_A)$. Since $\alpha_n$ and $x_A$ belong to the convex subcomplex $M_A$, we deduce that $\mu(\alpha_n,\beta_n,x_A) \in M_A$, so $d_1(m_n,M_A) \leq R$. For the same reason, we have $d_1(m_n,M_B) \leq R$. As a consequence, we have $\forall n \in \N, m_n \in S$.

Since $S$ is finite, up to passing to a subsequence we may assume that $\forall n \in \N, m_n=x_0$ is constant.

\mk

\begin{figure}[!h]
\def\svgwidth{4cm}
\center
\input{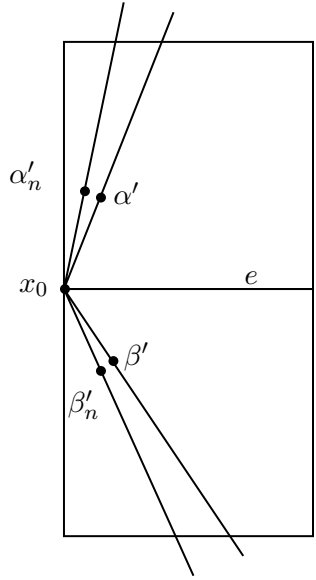}
\caption{The proof of Proposition~\ref{pro:2 abelian subgroups}}
\label{fig:cubes}
\end{figure}

Fix $\varepsilon>0$, and for each $n \in \N$, let $\alpha'_n$ (resp $\beta'_n$) be the point on the CAT(0) geodesic segment between $x_0$ and $\alpha_n$ (resp. $\beta_n$) at $d_2$ distance $\varepsilon$ from $x_0$ (see Figure~\ref{fig:cubes}).

We know that $\mu(x_0,\alpha'_n,\alpha_n)=\alpha'_n$, $\mu(x_0,\beta'_n,\beta_n)=\beta'_n$ and $\mu(x_0,\alpha_n,\beta_n)=x_0$. Hence by using Equation~(\ref{eqn:median}) from Section~\ref{subsec:median}, we deduce that
\beq
\mu(x_0,\alpha_n,\beta'_n) &=& \mu(x_0,\alpha_n,\mu(x_0,\beta_n,\beta'_n)) \\
&=& \mu(\mu(x_0,\alpha_n,x_0),\mu(x_0,\alpha_n,\beta_n),\beta'_n) = \mu(x_0,x_0,\beta'_n) = x_0, \mbox{ and }\\
\mu(x_0,\alpha'_n,\beta'_n) &=& \mu(x_0,\mu(x_0,\alpha'_n,\alpha_n),\beta'_n) \\
&=& \mu(\mu(x_0,x_0,\beta'_n),\mu(x_0,\alpha_n,\beta'_n),\alpha'_n) = \mu(x_0,x_0,\alpha'_n) = x_0.
\eeq

But the sequence $(\alpha'_n)_{n \in \N}$ (resp. $(\beta'_n)_{n \in \N}$) actually converges to the point $\alpha'$ (resp. $\beta'$) on the CAT(0) geodesic ray from $x_0$ to $a(\pif)$ (resp. $b(\pif)$) at $d_2$ distance $\varepsilon$ from $x_0$. Hence we conclude that $\mu(x_0,\alpha',\beta')=x_0$. In other words, the path $[\alpha',x_0] \cup [x_0,\beta']$ is monotone.

On the other hand, we have $\sphericalangle_{x_0}(\alpha',\beta') = \sphericalangle_{x_0}(a(\pif),b(\pif)) \leq \sphericalangle(a(\pif),b(\pif))$. By contradiction, assume that we have $\sphericalangle(a(\pif),b(\pif)) < \frac{\pi}{2}$, then $\sphericalangle_{x_0}(\alpha',\beta')<\f{\pi}{2}$. Let $\sigma_{\alpha'},\sigma_{\beta'}$ denote the minimal (closed) simplices in the link $L$ of $X$ at $x_0$ containing $\alpha'$ and $\beta'$ respectively. Since $\sphericalangle_{x_0}(\alpha',\beta')<\f{\pi}{2}$, we know that $d(\sigma_{\alpha'},\sigma_{\beta'}) < \f{\pi}{2}$.

\mk

We claim that $\sigma_{\alpha'}$ and $\sigma_{\beta'}$ are not disjoint. Consider any two disjoint simplices $\sigma,\sigma'$ in $L$, we will prove that $d(\sigma,\sigma') \geq \f{\pi}{2}$. By considering a geodesic from $\sigma$ to $\sigma'$ (if it exists, otherwise the distance is infinite), up to reducing the distance between $\sigma$ and $\sigma'$, we can assume that both $\sigma$ and $\sigma'$ are contained in a common simplex $\sigma_0$. We choose $\sigma_0$ minimal, so that $\sigma_0$ is isometric to the spherical join $\sigma \star \sigma'$. By definition of the spherical join, the distance between $\sigma$ and $\sigma'$ is equal to $\f{\pi}{2}$. So $\sigma_{\alpha'}$ and $\sigma_{\beta'}$ are not disjoint.

\mk

Therefore there exists an edge $e$ in $X$ adjacent to $x_0$ whose image in the link of $x_0$ belongs to $\sigma_{\alpha'} \cap \sigma_{\beta'}$. So $\sphericalangle_{x_0}(\alpha',e)<\f{\pi}{2}$ and $\sphericalangle_{x_0}(\beta',e)<\f{\pi}{2}$. If we consider a shifted hyperplane $H$ dual to $e$ close to $x_0$ (the CAT(0) orthogonal of $e$ at a point near $x_0$), we see that $H$ separates $x_0$ and $\{\alpha',\beta'\}$: this contradicts the monotonicity of the path $[\alpha',x_0] \cup [x_0,\beta']$.

As a consequence, we have $\sphericalangle(a(\pif),b(\pif)) \geq \frac{\pi}{2}$.
\ep


We can now prove the first two results giving obstructions to being virtually cocompactly cubulated. They show how to combine Propositions~\ref{pro:proper action of Im} and \ref{pro:2 abelian subgroups}.

Recall that a subgroup $H$ is said to \emph{virtually contain} $g$ if there exists $n \geq 1$ such that $g^n \in H$.

\blem \label{lem:non_existence} There is no group $G$ satisfying the following.
\bit
\item there exist elements $a,b \in G$ such that $\<a,b\> \simeq A(p)$, for some $p \geq 3$,
\item there exists a finite index normal subgroup $G_0$ of $G$ acting geometrically on a CAT(0) cube complex $X$,
\item there exists an abelian subgroup $C$ of $G_0$ virtually containing $z_{ab}$,
\item there exists an abelian subgroup $A$ of $G_0$ virtually containing $a$ such that $A \cap C$ is finite,
\item if $p$ is even, there exists an abelian subgroup $B$ of $G_0$ virtually containing $b$ such that $B \cap C$ is finite,
\item for every $g \in G$, the groups $gAg^{-1}$ and $gCg^{-1}$ (and $gBg^{-1}$ if $p$ is even) are convex-cocompact in $X$.
\eit
\elem

\bp By contradiction, assume that such a group $G$ exists. We will produce a proper action of $G$ on a CAT(0) space.

\mk

Consider the induced action of $G$ on the finite-dimensional CAT(0) cube complex $X^{G/G_0}$. To describe this action, one can for instance identify $X^{G/G_0}$ with the space of right $G_0$-equivariant maps from $G$ to $X$, endowed with the action of $G$ by left translations. This provides a proper action of $G$ on the CAT(0) cube complex $X^{G/G_0}$ by cubical isometries (this idea comes from~\cite[Remark~1]{bridson_mcg_cat0}).

\mk

Let $N \geq 1$ such that $a^N \in A$, and let $M \geq 1$ such that $z_{ab}^M \in C$. We will now prove that the attractive endpoints of $a^N$ and $z_{ab}^M$ in $\partial X^{G/G_0}$ satisfy $\sphericalangle_{\partial X^{G/G_0}}(a^N(\pif),z_{ab}^M(\pif)) \geq \f{\pi}{2}$, which will contradict Proposition~\ref{pro:proper action of Im}.

\mk

Let us denote $G/G_0=\{g_1G_0,\dots,g_nG_0\}$. The action of $G_0$ on $X^{G/G_0}$ preserves each factor, and the action of $G_0$ on $X^{g_iG_0}$ is isomorphic to the conjugate by $g_i$ of the original action of $G_0$ on $X$. For each $1 \leq i \leq n$, we know that $g_i a^N g_i^{-1} \in g_i A g_i^{-1}$, and the subgroup $g_i A g_i^{-1}$ is convex-cocompact in $X$ by assumption. Similarly, for each $1 \leq i \leq n$, we know that $g_i z_{ab}^M g_i^{-1} \in g_i C g_i^{-1}$, and the subgroup $g_i C g_i^{-1}$ is convex-cocompact in $X$ by assumption. Furthermore, the intersection $g_i A g_i^{-1} \cap g_i C g_i^{-1} = g_i (A \cap C) g_i^{-1} = \{1\}$ is finite. According to Proposition~\ref{pro:2 abelian subgroups}, we deduce that $\sphericalangle_{\partial X^{g_iG_0}}(a^N(\pif),z_{ab}^M(\pif)) \geq \f{\pi}{2}$.

\mk

Note that the visual boundary of the finite product $X^{G/G_0}$ is isometric to the spherical join of the visual boundaries $\partial X^{g_iG_0}$ of each factor (see~\cite[Proposition~I.5.15]{bridson_haefliger}). Furthermore, by definition of the distance on the spherical join $S_1 \ast S_2$ of two metric spaces $S_1,S_2$  (see~\cite[Definition~I.5.13]{bridson_haefliger}), if $x,y \in S_1 \ast S_2$ are such that their distances in each factor are at least $\frac{\pi}{2}$, then $x$ and $y$ are at distance at least $\frac{\pi}{2}$.

Since for any $1 \leq i \leq n$ we have $\sphericalangle_{\partial X^{g_iG_0}}(a^N(\pif),z_{ab}^M(\pif)) \geq \f{\pi}{2}$, we deduce that $\sphericalangle_{\partial X^{G/G_0}}(a^N(\pif),z_{ab}^M(\pif)) \geq \f{\pi}{2}$. By symmetry if $p$ is even, we also have $\sphericalangle_{\partial X^{G/G_0}}(b^N(\pif),z_{ab}^M(\pif)) \geq \f{\pi}{2}$.

\mk

Remark now that the group $G$ contains the dihedral Artin group $\<a,b\>$, and acts properly by semisimple isometries on the CAT(0) cube complex $X^{G/G_0}$, so this contradicts Proposition~\ref{pro:proper action of Im}. This concludes the proof.
\ep

\blem \label{lem:non_existence_step2} There is no group $G$ satisfying the following.
\bit
\item there exist elements $a,b \in G$ such that $\<a,b\> \simeq A(p)$, for some $p \geq 3$,
\item there exists a finite index normal subgroup $G_0$ of $G$ acting geometrically on a CAT(0) cube complex $X$,
\item there exists an abelian subgroup $C$ of $G_0$ commuting with $a$ and $b$ and virtually containing $z_{ab}$,
\item there exists an abelian subgroup $A$ of $G_0$ virtually containing $a$ such that $A \cap \<z_{ab}\>=\{1\}$,
\item if $p$ is even, there exists an abelian subgroup $B$ of $G_0$ virtually containing $b$ such that $AB \cap \<z_{ab}\>=\{1\}$,
\item for every $g \in G$, the groups $gAg^{-1}$ and $gCg^{-1}$ (and $gBg^{-1}$ if $p$ is even) are convex-cocompact in $X$.
\eit
\elem

\bp
By contradiction, assume that there exists a counterexample $G$. Assume furthermore that, among all counterexamples, the dimension of the CAT(0) cube complex $X$ is minimal. We will prove that $G$ is then a counterexample to Lemma~\ref{lem:non_existence}.

\mk

We will now prove that $A \cap C$ is finite. Let $W=A \cap C$. Since $W$ is abelian, according to Proposition~\ref{pro:cubulation_centralizer}, there exists a finite index subgroup $W_0 \subset W \cap G_0$ such that the centralizer $H=Z_G(W_0)$ has a finite index normal subgroup $H_0 \subset G_0$ acting geometrically on a CAT(0) cube complex $Y$ with $\dim Y \leq \dim X$. Furthermore, for every $g \in H$, the groups $gAg^{-1} \cap H_0$ and $gCg^{-1} \cap H_0$ (and $gBg^{-1} \cap H_0$ if $p$ is even) are convex-cocompact in $Y$. Also note that, since $a,b$ commute with $C$ and $W \subset C$, we deduce that $a,b \in H$.

\mk

According to Proposition~\ref{pro:cubulation_central_quotient}, the group $W_0$ has a finite index subgroup $W'$ such that $G'=H/W'$ has a finite index subgroup $G'_0 = H_0/W'$ acting geometrically on a CAT(0) cube complex $X'$, with $\dim X' \leq \dim Y - \rk W' \leq \dim X - \rk W'$. Furthermore, for every $g \in H$, the groups $(gAW'g^{-1} \cap H_0)/W'$ and $(gCW'g^{-1} \cap H_0)/W'$ (and $(gBW'g^{-1} \cap H_0)/W'$ if $p$ is even) are convex-cocompact in $X'$. 

\mk

We will prove that the group $G'$ is also a counterexample to Lemma~\ref{lem:non_existence_step2}.

\bit
\item We will prove that $a'=aW'$ and $b'=bW'$ span a subgroup of $G'$ isomorphic to $\<a,b\> \simeq A(p)$. Since $W'$ is central in $H$, we deduce that $\<a,b\> \cap W'$ is central in $\<a,b\>$, so $\<a,b\> \cap W' \subset \<z_{ab}\> \cap  A = \{1\}$. Hence $\<a,b\> \cap W' = \{1\}$, and $\<a',b'\> \simeq \<a,b\> \simeq A(p)$.
\item The finite index subgroup $G'_0$ of $G'$ acts geometrically on the CAT(0) cube complex $X'$.
\item The abelian subgroup $C'=(CW' \cap H)/W'$ of $G'$ virtually contains $z_{a'b'}$, since both $C$ and $H$ virtually contain $z_{ab}$. Furthermore, since $a$ and $b$ commute with $C$, we deduce that $a'$ and $b'$ commute with $C'$.
\item The abelian subgroup $A'=(AW' \cap H) / W'$ of $G'$ virtually contains $a'$, since both $A$ and $H$ virtually contain $a$. Furthermore, since $W' \subset A$, we have $AW' \cap (\<z_{ab}\>W') = A \cap (\<z_{ab}\>W') = (A \cap \<z_{ab}\>)W'=W'$, so $A' \cap \<z_{a'b'}\>=\{1\}$.
\item If $p$ is even, similarly the abelian subgroup $B'=(BW' \cap H) / W'$ of $G'$ virtually contains $b'$. Furthermore, assume that $n \in \Z$ and $w \in W'$ are such that $z_{a'b'}^nw \in (AW' \cap H)(BW' \cap H)$, we will prove that $n=0$. As $W' \subset BW' \cap H$, we have $z_{a'b'}^n \in (AW' \cap H)(BW' \cap H)$. Since $W' \subset A$, we have $AW'=A$ so $z_{a'b'}^n \in (A \cap H)(BW' \cap H)$. As $W' \subset BW' \cap H$ and $W'$ commutes with $H$, we have $z_{a'b'}^n \in (A \cap H)W'(B \cap H)$. Since $W' \subset A \cap H$, we have $z_{a'b'}^n \in (A \cap H)(B \cap H) \subset AB$. By assumption, $n=0$. So $A'B' \cap \<z_{a'b'}\>=\{1\}$.
\item For any $gW' \in G'$, we have seen that the groups $(gW') A'(gW')^{-1} \cap G'_0= (gAW'g^{-1} \cap H_0)/W'$ and $(gW')C'(gW')^{-1}  \cap G'_0= (gCW'g^{-1} \cap H_0)/W'$ (and  $(gW') B'(gW')^{-1}  \cap G'_0= (gBW'g^{-1} \cap H_0)/W$ if $p$ is even) are convex-cocompact in $X'$. 
\eit

As a consequence, the group $G'$ is also a counterexample to Lemma~\ref{lem:non_existence_step2}, with $\dim X' \leq \dim X - \rk W'$. By minimality of $\dim X$, we deduce that $\rk W'=0$, so $W=A \cap C$ is finite.

\mk

If $p$ is even, we argue similarly that $B \cap C$ is finite. To be precise, the only difference in the statement between $A$ and $B$ is that $A \cap \<z_{ab}\> =\{1\}$ and $AB \cap \<z_{ab}\> = \{1\}$. But this assumption is, in fact, symmetrical with respect to $A$ and $B$: indeed $AB \cap \<z_{ab}\> = \{1\}$ implies both $B \cap \<z_{ab}\>=\{1\}$ and $BA \cap \<z_{ab}\> = \{1\}$.

\mk

As a consequence, the group $G$ contradicts Lemma~\ref{lem:non_existence_step2}. Therefore there exists no such group $G$.
\ep

We are now ready to state the most general and self-contained result giving an obstruction to being virtually cocompactly cubulated. In the proof, we will produce small convex-cocompact subgroups using highest abelian subgroups.

\bpro \label{pro:not_cubulated} Let $G$ be a group satisfying the following.
\bit
\item There exist elements $a,b \in G$ such that $\<a,b\> \simeq A(p)$, for some $p \geq 3$,
\item For every $n \geq 1$, we have $Z_G(a^n)=Z_G(a)$ and $Z_G(b^n)=Z_G(b)$,
\item There exists $\alpha \in G$ commuting with $a$, such that no non-zero powers of $\alpha$ and $z_{ab}$ commute,
\item If $p$ is even, there exists $\beta \in G$, commuting with $b$, such that for every $q \geq 1$ we have $Z_G(a,z_{ab}^q,\alpha^q) Z_G(b,\beta^q) \cap \<z_{ab}\> = \{1\}$.
\eit
Then $G$ is not virtually cocompactly cubulated.
\epro

\bp
By contradiction, assume that some finite index normal subgroup $G_0$ of $G$ acts geometrically on a CAT(0) cube complex $X$. We will prove that $G$ is then a counterexample to Lemma~\ref{lem:non_existence_step2}.

\mk

Let $A$ denote the intersection of two highest maximal abelian subgroups of $G$ virtually containing $\{a,\alpha\}$ and $\{a,z_{ab}\}$ respectively. Since for every $n \geq 1$ we have $Z_G(a^n)=Z_G(a)$, by maximality we deduce that $a \in A$. Similarly, let $B$ denote a highest maximal abelian subgroup of $G$ virtually containing $b$ and $\beta$: $B$ contains $b$.

\mk

Let $A',B'$ denote two highest maximal abelian subgroups of $G$ virtually containing $\{a,z_{ab}\}$ and $\{b,z_{ab}\}$ respectively. As for $A$ and $B$, we know that $a \in A'$ and $b \in B'$. So $C=A' \cap B'$ commutes with $a$ and $b$. Furthermore $C$ virtually contains $z_{ab}$.

\mk

According to Theorem~\ref{thm:flat_torus} and Lemma~\ref{lem:intersection of convex subgroups}, we know that for every $g \in G$, the groups $gAg^{-1} \cap G_0$ and $gCg^{-1} \cap G_0$ (and $gBg^{-1} \cap G_0$ if $p$ is even) are convex-cocompact in $X$.

\mk

By definition of $A$, there exists $m \geq 1$ such that $\alpha^m$ commutes with $A$. Since no non-zero powers of $\alpha$ and $z_{ab}$ commute, we deduce that $A \cap \<z_{ab}\> = \{1\}$.

\mk

If $p$ is even, then by definition of $A$ and $B$ we know that there exists $q \geq 1$ such that $A \subset Z_G(a,z_{ab}^q,\alpha^q)$ and $B \subset Z_G(b,\beta^q)$, hence $\<z_{ab}\> \cap AB = \{1\}$.

\mk

In conclusion, $G$ is a counterexample to Lemma~\ref{lem:non_existence_step2}. Therefore $G$ is not virtually cocompactly cubulated.
\ep

In order to apply Proposition~\ref{pro:not_cubulated} to Artin groups, we need the following technical results.

\blem \label{lem:abc_odd}
Assume that $S=\{a,b,c\}$, and $M=(m_{st})_{s,t \in S}$ is a Coxeter matrix such that $m_{ab}$ is finite and odd, $m_{ac}$ is finite and $m_{bc}$ is different from $2$. Then $z_{ac}$ and $z_{ab}$ do not virtually commute.
\elem

\bp
Assume there exist $n,m \in \Z$ such that $z_{ab}^n$ and $z_{ac}^m$ commute. We can assume that $n,m \geq 0$. Then $z_{ab}^nz_{ac}^m=z_{ac}^mz_{ab}^n$ is an equality between positive words, so by~\cite{paris_monoid} they are equal in the positive monoid: one can pass from one to the other by applying only the standard relations of $A(S)$. But the relation between $b$ and $c$ cannot be used since $m_{bc} \neq 2$, the subword $w_{m_{bc}}(b,c)$ cannot appear. As a consequence, starting from $z_{ab}^nz_{ac}^m$ it is not possible to obtain a word with a letter $c$ on the left of a letter $b$. This implies that $n=0$ or $m=0$: no non-trivial powers of $z_{ab}$ and $z_{ac}$ commute.
\ep

\blem \label{lem:abc}
Assume that $S=\{a,b,c\}$, and $M=(m_{st})_{s,t \in S}$ is a Coxeter matrix with finite entries such that $m_{ab}$ and $m_{ac}$ are even numbers different from $2$, and $m_{bc}$ is even. Then $z_{ac}$ and $z_{ab}$ do not virtually commute. 

Furthermore, fix $q \geq 1$ and assume that $u \in A(M)$ commutes with $a$, $z_{ab}^q$ and $z_{ac}^q$, and that $v \in A(M)$ commutes with $b$ and $z_{bc}^q$. If there exists $n \in \Z$ such that $uv=z_{ab}^n$, then $n=0$.
\elem

\bp
Assume there exist $n,m \in \Z$ such that $z_{ab}^n$ and $z_{ac}^m$ commute. We can assume that $n,m \geq 0$. Then $z_{ab}^nz_{ac}^m=z_{ac}^mz_{ab}^n$ is an equality between positive words, so by~\cite{paris_monoid} they are equal in the positive monoid: one can pass from one to the other by applying only the standard relations of $A(S)$. But, starting from $z_{ab}^nz_{ac}^m$, the letters $c,a,b$ cannot appear in that order, since $m_{ab} \geq 4$ and $m_{ac} \geq 4$. This implies that $n=0$ or $m=0$: no non-trivial powers of $z_{ab}$ and $z_{ac}$ commute.

\mk

Since $\f{1}{m_{ab}}+\f{1}{m_{bc}}+\f{1}{m_{ac}} \leq 1$, the group $A(M)$ is not of spherical type. According to~\cite[Theorem~B]{charney_davis_kpi1}, the Deligne complex, with Moussong's metric, is CAT(0). According to~\cite[Theorem~1]{godelle_cat0}, Godelle's Property $(\star \star \star)$ can be applied to the pair $(\{a\},\{a\})$ to conclude that the centralizer of $a$ in $A(M)$ can be described, using ribbons, as
$$Z_{A(M)}(a) = \<a,z_{ab},z_{ac}\>.$$ 

\mk

By Charney and Davis (see~\cite[Theorem~B]{charney_davis_kpi1} and \cite[Corollary~1.4.2]{charney_davis_salvetti}), the cohomological dimension of $A(M)$ is $2$, so the maximal rank of an abelian subgroup of $A(M)$ is $2$. As a consequence, the only elements of $Z_{A(M)}(a)$ commuting with $z_{ab}^q$ and $z_{ac}^q$ are powers of $a$, so $u \in \<a\>$. Let $r \in \Z$ such that $u=a^r$.

\mk

Assume that $v \in A(M)$ commutes with $b$ and $z_{bc}^q$, and that $n \in \Z$ is such that $uv=z_{ab}^n$. We will prove that $n=0$. Remark that $z_{bc}^q$ commutes with $v=u^{-1}z_{ab}^n=a^{-r}z_{ab}^n$.

Consider the homomorphism $\phi_c: A(M) \ra \<a,b\>$ sending $a,b$ to $a,b$ and sending $c$ to $1$. Since all integers defining $M$ are even, $\phi_c$ is a well-defined group homomorphism. Then $\phi_c(z_{bc}^q) = b^{\frac{qm_{bc}}{2}}$ commutes with $\phi_c(a^{-r}z_{ab}^n)=a^{-r}z_{ab}^n$, so $r=0$.

We have $u=a^r=1$, so $z_{ab}^n=v$ commutes with $z_{bc}^q$, which implies that $n=0$.
\ep

\blem \label{lem:abcd}
Assume that $S=\{a,b,c,d\}$, and $M=(m_{st})_{s,t \in S}$ is a Coxeter matrix with all entries even or infinite. Assume furthermore that $m_{ab}$ is finite and different from $2$, $m_{ac}$ and $m_{bd}$ are finite, and $m_{ad},m_{bc}$ are different from $2$. Then $z_{ac}$ and $z_{ab}$ do not virtually commute. 

Furthermore, fix $q \geq 1$ and assume that $u \in A(M)$ commutes with $a$, $z_{ab}^q$ and $z_{ac}^q$, and that $v \in A(M)$ commutes with $b$ and $z_{bd}^q$. If there exists $n \in \Z$ such that $uv=z_{ab}^n$, then $n=0$.\elem

\bp
Following the same proof as in Lemma~\ref{lem:abc_odd}, we see that since $m_{bc} \neq 2$, no non-zero powers of $z_{ab}$ and $z_{ac}$ commute.

\mk

Without loss of generality, we can assume up to passing to the corresponding quotient that $m_{ac}=m_{bd}=m_{cd}=2$. Since all finite entries of the Coxeter matrix are even, every irreducible spherical parabolic subgroup of $A(M)$ has rank $1$ or $2$. According to~\cite[Theorem~B]{charney_davis_kpi1}, the Deligne complex, with Moussong's metric, is CAT(0). According to~\cite[Theorem~1]{godelle_cat0}, Godelle's Properties $(\star)$ and $(\star \star \star)$ can be applied to the pair $(\{a,c\},\{a,c\})$ to conclude that the commensurator of $\<a,c\>$ in $A(M)$ is equal to $\<a,c\>$.

Since $u$ commutes with $a$ and $c^q$, $u$ commensurates $\<a,c\>$, so $u \in \<a,c\>$. As $m_{bc} \neq 2$, no non-zero powers of $c$ and $z_{ab}$ commute, so $u \in \<a\>$. Let $r \in \Z$ such that $u=a^r$.

\mk

Assume that $v \in A(M)$ commutes with $b$ and $z_{bd}^q$, and that $n \in \Z$ is such that $uv=z_{ab}^n$. We will prove that $n=0$. Remark that $z_{bd}^q$ commutes with $v=u^{-1}z_{ab}^n=a^{-r}z_{ab}^n$.

Consider the homomorphism $\phi_d: A(M) \ra \<a,b,c\>$ sending $a,b,c$ to $a,b,c$ and sending $d$ to $1$. We deduce that $\phi_d(z_{bd}^q) = \phi_d(b^qd^q)=b^q$ commutes with $\phi_d(a^{-r}z_{ab}^n)=a^{-r}z_{ab}^n$, so $r=0$.

We have $u=a^r=1$, so $z_{ab}^n=v$ commutes with $z_{bd}^q$, which implies that $n=0$.
\ep

We will now prove that the statements of Conjecture~\ref{conj:main} and Conjecture~\ref{conj:main2} are equivalent.

\bpro
Let $M=(m_{st})_{s,t \in S}$ be a finite Coxeter matrix. Consider the following five conditions.
\begin{enumerate}[label=\Alph*.]
\item For each pairwise distinct $a,b,c \in S$ such that $m_{ab}$ is odd, either $m_{ac}=m_{bc}=\infty$ or $m_{ac}=m_{bc}=2$.
\item For each distinct $a,b \in S$ such that $m_{ab}$ is even and different from $2$, there is an ordering of $\{a,b\}$ (say $a < b$) such that, for every $c \in S \bs \{a,b\}$, one of the following holds:
\ben
\item $m_{ac}=m_{bc}=2$,
\item $m_{ac}=2$ and $m_{bc}=\infty$,
\item $m_{ac}=m_{bc}=\infty$, or
\item $m_{ac}$ is even and different from $2$, $a<c$ in the ordering of $\{a,c\}$, and $m_{bc}=\infty$.
\een
\een

\ben
\item There exist $3$ pairwise distinct $a,b,c \in S$ such that $m_{ab}$ is odd, $m_{ac} \neq \infty$ and $m_{bc} \neq 2$.
\item There exist $3$ pairwise distinct $a,b,c \in S$ such that $m_{ab}$ and $m_{ac}$ are even numbers different from $2$, and $m_{bc} \neq \infty$.
\item There exist $4$ pairwise distinct $a,b,c,d \in S$ such that $m_{ab} \not\in \{2,\infty\}$, $m_{ac}, m_{bd} \neq \infty$ and $m_{ad},m_{bc} \neq 2$.
\een
Then $A.$ and $B.$ hold if and only $1.$, $2.$ and $3.$ do not hold.
\epro

\bp
Assume first that $1.$, $2.$ or $3.$ holds, we will prove that $A.$ or $B.$ do not hold.
\ben
\item Assume that there exist $3$ pairwise distinct $a,b,c \in S$ such that $m_{ab}$ is odd, $m_{ac} \neq \infty$ and $m_{bc} \neq 2$. Then $a,b,c$ contradict Condition $A$.
\item Assume that there exist $3$ pairwise distinct $a,b,c \in S$ such that $m_{ab}$ and $m_{ac}$ are even numbers different from $2$, and $m_{bc} \neq \infty$. Then $a,b,c$ contradict Condition $B.(d)$.
\item Assume that there exist $4$ pairwise distinct $a,b,c,d \in S$ such that $m_{ab} \not\in \{2,\infty\}$ is even, $m_{ac}, m_{bd} \neq \infty$ and $m_{ad},m_{bc} \neq 2$. If an ordering of $\{a,b\}$ as in Condition $B.$ existed, we should have both $a<b$ and $b<a$, which is a contradiction.
\een

\mk

Assume now that $1.$, $2.$ and $3.$ do not hold, we will prove that $A.$ and $B.$ hold.
\ben[label=\Alph*.]
\item Consider three pairwise distinct $a,b,c \in S$ such that $m_{ab}$ is odd. Since Condition $1.$ does not hold, we have $m_{ac}=m_{bc} \in \{2,\infty\}$.
\item Consider distinct $a,b \in S$ such that $m_{ab}$ is even and different from $2$. If there exists $c \in S \bs \{a,b\}$ such that $m_{ac} \neq \infty$ and $m_{bc} = \infty$, choose the ordering $a<b$. If there exists $d \in S \bs \{a,b\}$ such that $m_{bd} \neq \infty$ and $m_{ad} = \infty$, choose the ordering $b<a$.If there is no such $c$ or $d$, choose an arbitrary ordering of $\{a,b\}$.

\mk

Notice that it is not possible that both $c$ and $d$ exist. By contradiction, assume that there exist $c,d \in S \bs \{a,b\}$ such that $m_{ac},m_{bd} \neq \infty$ and $m_{bc},m_{ad} = \infty$. This contradicts Condition $3$.

\mk

Now that the ordering of $\{a,b\}$ is well-defined, say $a<b$, we will check that it satisfies the required properties. Fix any $c \in S \bs \{a,b\}$.

\mk

Assume first that $m_{ac},m_{bc} \in \{2,\infty\}$. Then since $a<b$, we do not have both $m_{ac}=\infty$ and $m_{bc}=2$.

Assume now that $m_{bc} \not\in \{2,\infty\}$. Since Condition $1.$ does not hold, $m_{bc}$ is even. Since Condition $2.$ does not hold, we have $m_{ac}=\infty$. This contradicts $a<b$.

Assume finally that $m_{ac} \not\in \{2,\infty\}$. Since Condition $1.$ does not hold, $m_{ac}$ is even. Since Condition $2.$ does not hold, $m_{bc}=\infty$. This implies that $a<c$ in the ordering of $\{a,c\}$.

\mk

As a consequence, Conditions $A.$ and $B.$ are satisfied.
\een
\ep

Let us recall the definition of the property $(\dagger)$ needed to prove Conjecture~\ref{conj:main}. Let $M=(m_{ab})_{a,b \in S}$ be a finite Coxeter matrix. We say that the Artin-Tits group $A(M)$ satisfies property $(\dagger)$ if
$$\forall s \in S, \forall n \geq 1, Z_{A(M)}(s^n) = Z_{A(M)}(s).$$

We can now prove the following, which is a restatement of Theorem~\ref{thm:main2}.

\bthm \label{thm:artin_obstruction}
Let $M=(m_{ab})_{a,b \in S}$ be a finite Coxeter matrix such that the Artin-Tits group $A(M)$ satisfies property $(\dagger)$. Assume that at least one of the following holds
\bit
\item there exist $3$ pairwise distinct $a,b,c \in S$ such that $m_{ab}$ is odd, $m_{ac} \neq \infty$ and $m_{bc} \neq 2$,
\item there exist $3$ pairwise distinct $a,b,c \in S$ such that $m_{ab}$ and $m_{ac}$ are even numbers different from $2$, and $m_{bc} \neq \infty$, or
\item there exist $4$ pairwise distinct $a,b,c,d \in S$ such that $m_{ab} \not\in \{2,\infty\}$, $m_{ac}, m_{bd} \neq \infty$ and $m_{ad},m_{bc} \neq 2$.
\eit
Then the Artin-Tits group $A(M)$ is not virtually cocompactly cubulated.
\ethm

\bp
\bit
\item Assume first that there exist $3$ pairwise distinct $a,b,c \in S$ such that $m_{ab}$ is odd, $m_{ac} \neq \infty$ and $m_{bc} \neq 2$. Then by Lemma~\ref{lem:abc_odd}, the element $z_{ac}$ commutes with $a$, but $z_{ac}$ and $z_{ab}$ do not virtually commute. By Proposition~\ref{pro:not_cubulated} applied with $\alpha=z_{ac}$, $A(M)$ is not virtually cocompactly cubulated.

\mk

Assume now that this first situation does not occur. For the two remaining cases, we will apply the same strategy.

\item Assume that there exist $3$ pairwise distinct $a,b,c \in S$ such that $m_{ab}$ and $m_{ac}$ are even numbers different from $2$, and $m_{bc}$ is finite and even. Let $\alpha=z_{ac}$ and $\beta=z_{bc}$. According to Lemma~\ref{lem:abc}, $\alpha$ and $a$ commute, but no non-zero powers of $\alpha$ and $z_{ab}$ commute.

\mk

Fix $q \geq 1$, and assume that $u \in Z_{A(M)}(a,z_{ab}^q,\alpha^q)$, $v \in z_{A(M)}(b,\beta^q)$ and $n \in \Z$ are such that $uv=z_{ab}^n$. We will prove that $n=0$.

Consider the homomorphism $\phi: A(M) \ra \<a,b,c\>$ sending $a,b,c$ to $a,b,c$ and sending every element in $S \bs \{a,b,c\}$ to $1$. According to the first case of the proof, we can assume that for every $d \in S \bs \{a,b,c\}$, the exponents $m_{ad},m_{bd}$ and $m_{cd}$ are even or infinite. Hence $\phi$ is a well-defined group homomorphism. Hence $\phi(u) \phi(v) = \phi(z_{ab}^n)=z_{ab}^n$.

Then $\phi(u)$ commutes with $\phi(a)=a$, $\phi(z_{ab}^q)=z_{ab}^q$ and $\phi(\alpha^q)=\phi(z_{ac}^q)=z_{ac}^q$, and $\phi(v)$ commutes with $\phi(b)=b$ and $\phi(\beta^q)=\phi(z_{bc}^q)=z_{bc}^q$. According to Lemma~\ref{lem:abc} applied to $\phi(u)$ and $\phi(v)$, we have $n=0$.

\mk

According to Proposition~\ref{pro:not_cubulated}, the group $A(M)$ is not virtually cocompactly cubulated.

\item Assume that there exist $4$ pairwise distinct $a,b,c,d \in S$ such that $m_{ab} \not\in \{2,\infty\}$, $m_{ac}, m_{bd} \neq \infty$ and $m_{ad},m_{bc} \neq 2$. Let $\alpha=z_{ac}$ and $\beta=z_{bd}$. According to Lemma~\ref{lem:abcd}, $\alpha$ and $a$ commute, but no non-zero powers of $\alpha$ and $z_{ab}$ commute.

\mk

Fix $q \geq 1$, and assume that $u \in Z_{A(M)}(a,z_{ab}^q,\alpha^q)$, $v \in z_{A(M)}(b,\beta^q)$ and $n \in \Z$ are such that $uv=z_{ab}^n$. We will prove that $n=0$.

Consider the homomorphism $\phi: A(M) \ra \<a,b,c,d\>$ sending $a,b,c,d$ to $a,b,c,d$ and sending every element in $S \bs \{a,b,c,d\}$ to $1$. According to the first case of the proof, we can assume that for every $t \in S \bs \{a,b,c,d\}$, the exponents $m_{at},m_{bt},m_{ct}$ and $m_{dt}$ are even or infinite. Hence $\phi$ is a well-defined group homomorphism. Hence $\phi(u) \phi(v) = \phi(z_{ab}^n)=z_{ab}^n$.

Then $\phi(u)$ commutes with $\phi(a)=a$, $\phi(z_{ab}^q)=z_{ab}^q$ and $\phi(\alpha^q)=\phi(z_{ac}^q)=z_{ac}^q$, and $\phi(v)$ commutes with $\phi(b)=b$ and $\phi(\beta^q)=\phi(z_{bd}^q)=z_{bd}^q$. According to Lemma~\ref{lem:abcd} applied to $\phi(u)$ and $\phi(v)$, we have $n=0$.

\mk

According to Proposition~\ref{pro:not_cubulated}, the group $A(M)$ is not virtually cocompactly cubulated.
\eit
\ep

\section{Cubulation of Artin groups}

\label{sec:cubulation}

\subsection{Cubulation of dihedral Artin groups}

Brady and McCammond showed (see~\cite{brady_mccammond_artinthree}) that for all $p \in \{2,\ldots,\infty\}$, the dihedral Artin group $A(p)$ is cocompactly cubulated. Let us recall their construction, which will be useful. We will need this construction when $p \not\in \{2,\infty\}$, but it works as well when $p=2$, so let us fix $p \neq \infty$ (when $p=\infty$, the Artin group is just the rank $2$ free group).

The Artin group $A(p)$ has the following presentation, due to Brady and McCammond:
$$A(p) = \< x,a_1,\ldots,a_p \st \forall 1 \leq i \leq p, a_ia_{i+1}=x\>,$$
where $a_{p+1}=a_1$. This can easily be seen, with $a_1$ and $a_2$ corresponding to the standard generators of $A(p)$. The presentation $2$-complex $K$ is a $K(\pi,1)$ for $A(p)$ consisting of $1$ vertex $v$, $p+1$ loops $a_1,\ldots,a_p,x$ and $p$ triangles $a_1a_2x^{-1},\ldots,a_pa_1x^{-1}$ (see Figure~\ref{fig:K}). 

\begin{figure}[!h]
\def\svgwidth{12cm}
\center
\input{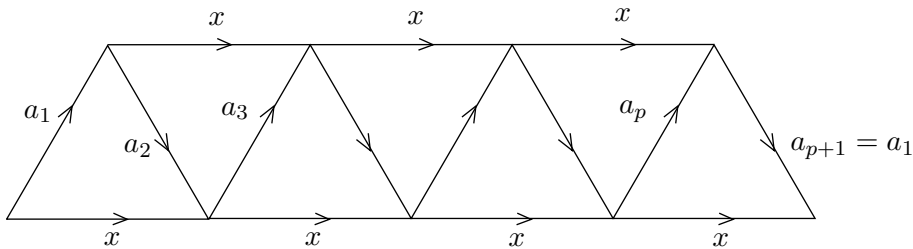}
\caption{Brady and McCammond's presentation $2$-complex $K$}
\label{fig:K}
\end{figure}

We will define another $K(\pi,1)$ for $A(p)$, which will be cubical and will have the same underlying topological space as $K$. Start with two vertices $v$ and $w$, and $p+2$ oriented edges between $v$ and $w$ labelled $\alpha_1,\ldots,\alpha_p,\beta^{-1},\gamma$. Finally, add the $p$ squares with boundary labeled by $\alpha_1\beta^{-1}\alpha_2\gamma^{-1},\ldots,\alpha_p\beta^{-1}\alpha_1\gamma^{-1}$ and let $X(A(p))$ denote the resulting cube complex. It is easy to see that the underlying topological space of $X(A(p))$ is homoeomorphic to $K$: $w$ corresponds to the midpoint of the edge $x$, the edge $x$ corresponds to the path $\gamma\beta^{-1}$, and each square corresponds to the union of the halves of two triangles of $K$ (see Figure~\ref{fig:Kbis}).

\begin{figure}[!h]
\def\svgwidth{12cm}
\center
\input{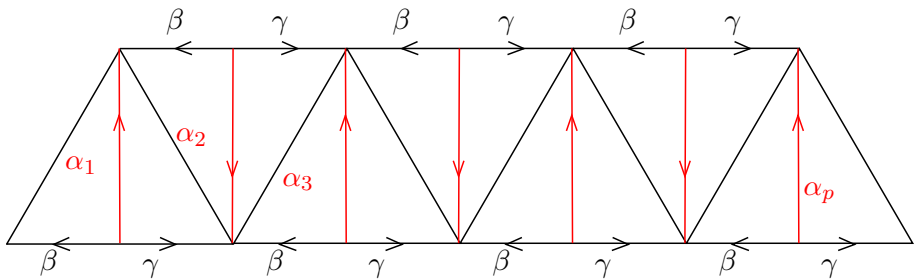}
\caption{The square complex $X(A(p))$}
\label{fig:Kbis}
\end{figure}

Hence $X(A(p))$ is also a $K(\pi,1)$ for $A(p)$. Furthermore, one has a complete description of the link of the vertex $v$: it has $p+2$ vertices labeled $\alpha_1,\ldots,\alpha_p,\beta^{-1},\gamma$, and it is the complete bipartite graph on $\{\alpha_1,\dots,\alpha_p\}$ and $\{\beta^{-1},\gamma\}$. This graph has no double edges and no triangle, so it is a flag simplicial complex. The description of the link of $w$ is similar, with all orientations reversed, so we deduce that $X(A(p))$ is a locally CAT(0) square complex. In particular, $A(p)$ is cocompactly cubulated.

\brk Notice that $X(A(p))$ is naturally isometric to the product of $\R$ and the infinite $p$-regular tree. In the case of the $3$-strand braid group $B_3\simeq A(3)$, one recovers in the central quotient the action of $B_3/Z(B_3) \simeq \PSL(2,\Z) \simeq \Z/2\Z * \Z/3\Z$ on its Bass-Serre $3$-regular tree. \erk

\subsection{Recubulation of even dihedral Artin groups}

\label{subsec:cubulation_even}

In the case where $p$ is even, there are two other natural CAT(0) square complexes on which the dihedral Artin group $A(p)$ acts geometrically. Each will be associated with one of the two generators $a$,$b$ of $A(p)$. We will describe the first one, associated with $a=a_1$.

\mk

Start with the same presentation $2$-complex $K$ as before, and remove all edges $a_2,a_4,\dots,a_p$ with even labels, and replace each pair of triangles $(a_{2i+1}a_{2i+2}x^{-1},a_{2i+2}a_{2i+3}x^{-1})$, for $0 \leq i \leq \f{p}{2}-1$, by a square with edges $a_{2i+1}xa_{2i+3}^{-1}x^{-1}$. We obtain a square complex $X_a(A(p))$ with one vertex $v$, $\f{p}{2}+1$ edges $x,a_1,a_3,\dots,a_{p-1}$ and $\f{p}{2}$ squares 
$a_1x{a_3}^{-1}x^{-1},\dots,a_{p-1}x{a_1}^{-1}x^{-1}$ (see Figure~\ref{fig:Keven}).

\begin{figure}[!h]
\def\svgwidth{9cm}
\center
\input{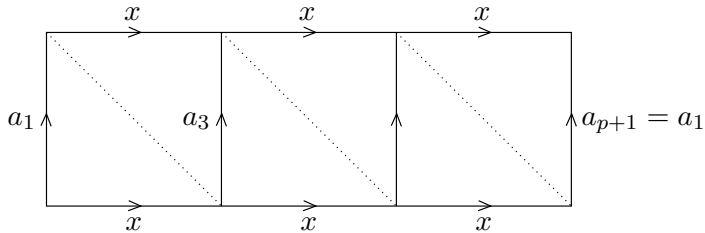}
\caption{The square complex $X_a(A(p))$}
\label{fig:Keven}
\end{figure}

The underlying topological space of $X_a(A(p))$ is $K$, so it is also a $K(\pi,1)$ for $A(p)$. Furthermore, one has a complete description of the link of the vertex $v$: it has $p+2$ vertices labeled $x,x^{-1},a_1,a_1^{-1},a_3,a_3^{-1},\dots,a_{p-1},a_{p-1}^{-1}$, and it is the complete bipartite graph on $\{x,x^{-1}\}$ and $\{a_1,a_1^{-1},a_3,a_3^{-1},\dots,a_{p-1},a_{p-1}^{-1}\}$. This graph has no double edges and no triangle, so it is a flag simplicial complex. So we deduce that $X_a(A(p))$ is a locally CAT(0) square complex.

\mk

The other locally CAT(0) square complex, denoted $X_b(A(p))$, is obtained by keeping only the edges with even labels and removing those with odd labels.

\mk

The fundamental difference of $X_a(A(p))$ and $X(A(p))$ is that, in the universal covers, the visual angles between the attractive fixed points of $a$ and $z_{ab}$ differ: in $X(A(p))$ that angle is acute, while in $X_a(A(p))$ it is equal to $\f{\pi}{2}$. This is due to the fact that, in $X_a(A(p))$, the edge $a=a_1$ belongs to the complex, so the subgroup $\<a\>$ is convex-cocompact in $X_a(A(p))$ but not in $X(A(p))$. This illustrates the case where $p$ is even in Proposition~\ref{pro:proper action of Im}.

\mk

In view of Lemma~\ref{lem:non_existence}, one can see that it is not possible to find a CAT(0) cube complex with a geometric action of $A(p)$ where both $\<a\>$ and $\<b\>$ are convex-cocompact. And if $p$ is odd, it is not even possible to find one where either $\<a\>$ or $\<b\>$ is convex-cocompact.

\subsection{Cubulation of Artin groups of even stars}

Let $M=(m_{ab})_{a,b \in S}$ be a finite Coxeter matrix, which is an ``even star'': there exists a ``central vertex'' $a \in S$ such that $\forall b,c \in S \bs \{a\}, m_{bc}=\infty$ and $\forall b \in S\bs \{a\}, m_{ab}$ is even.

\mk

We will now prove a particular case of Theorem~\ref{thm:main1}, namely showing that $A(M)$ is cocompactly cubulated. Note that J.~Huang, K.~Jankiewicz and P.~Przytycki independently gave the same construction in~\cite{huang_jankiewicz_przytycki}.

\mk

Write $S \bs \{a\}=\{b_1,\dots,b_m\}$. For each $1 \leq i \leq m$, the subgroup $A(\{ab_i\})$ of $A(M)$ spanned by $a$ and $b_i$ is a dihedral Artin group with even integer: let $X_a(A(\{ab_i\}))$ denote the previously constructed locally CAT(0) square complex with fundamental group $A(\{ab_i\})$, where some edge $e_i$ in $X_a(A(\{ab_i\}))$ represents $a$.

\mk

Consider now the square complex $X(A(M))$ which is the gluing of the square complexes $X_a(A(\{ab_1\})), \dots ,X_a(A(\{ab_m\}))$ where all edges $e_1,\dots,e_m$ are identified with a single edge $e$. By Van Kampen Theorem, the fundamental group of $X(A(M))$ is the free product of $A(\{ab_1\}), \dots, A(\{ab_m\})$ amalgamated over the cyclic subgroup $\<a\>$, which is precisely isomorphic to the Artin group $A(M)$.

\mk

If two squares $Q,Q'$ in $X(A(M))$ which do not lie in the same $X_a(A(\{ab_i\}))$ share an edge, this edge is $e$. As a consequence, for any triple of square $Q,Q',Q''$ of $X(A(M))$ which do not lie in the same $X_a(A(\{ab_i\}))$, if they pairwise share an edge, then their triple intersection is $e$. So $X(A(M))$ is a locally CAT(0) square complex. As a consequence, $A(M)$ is cocompactly cubulated.

\subsection{General case}

Let $M=(m_{ab})_{a,b \in S}$ be a finite Coxeter matrix satisfying the two conditions of Conjecture~\ref{conj:main}. We will show that the Artin group $A(M)$ is cocompactly cubulated.

Let $S_0 = \{a \in S \st \forall b \in S \bs \{a\}, m_{ab} \in \{2,\infty\}\}$ denote the set of vertices having all incident labels equal to $2$ or $\infty$.

Let $S_1=\{a_1,b_1\},\ldots,S_n=\{a_n,b_n\}$ denote the pairs of vertices of $S$ for which the edge $a_ib_i$ has an odd label, for $1 \leq i \leq n$ (possibly $n=0$).

Let $S_{n+1},\dots,S_{n+p}$ denote the subsets of vertices of $S \bs S_0$ for which the induced matrix $M_{S_{n+i} \times S_{n+i}}$ is an even star with central vertex $a_i \in S_{n+i}$, for $1 \leq i \leq p$ (possibly $p=0$). 

By assumption, we have $S=\bigsqcup_{0 \leq i \leq n+p} S_i$. 

\bigskip

We will consider cube complexes with edges labeled in ${\cal P}(S)$, the power set of $S$.

\mk

Let $X_0$ be the Salvetti cube complex of the right-angled Artin group of the graph induced by $S_0$: we will recall here its construction (see~\cite{salvetti}). It has one vertex and its edge set is $S_0$: each edge is labeled by $\{a\}$, for some $a \in S_0$. For each simplex $T \subset S_0$, we add a $|T|$-cube, by identifying each of the $|T|$ parallel classes of edges of $[0,1]^{|T|}$ with the edges labeled by $\{t\}$, for $t \in T$. Then by Theorem~\ref{thm:gromov criterion}, $X_0$ is locally CAT(0) cube complex.

For each $1 \leq i \leq n$, let $X_i$ denote a copy of the previously constructed cube complex $X(A(p_i))$ for the subgroup generated by $a_i$ and $b_i$, where $p_i$ is odd. Label each edge of $X_i$ by $\{a_i,b_i\}$.

For each $n+1 \leq i \leq n+p$, let $X_i$ denote a copy of the previously constructed cube complex $X(A(S_i))$ for the subgroup generated by $S_i$. Label the edge corresponding to the element $a_i$ by $\{a_i\}$, and label each other edge coming from the square complex $X_{a_i}(A(\{a_ib\}))$ by $\{a_i,b\}$, for every $b \in S_i \bs \{a_i\}$.

\mk

Consider the following cube complex $X$, which will be a cube subcomplex of the direct product $\prod_{i=0}^{n+p} X_i$. For each set of cubes $Q_0,\dots,Q_{n+p}$ of $X_0,\dots,X_{n+p}$ respectively, we will add the cube $Q_0 \times \dots \times Q_{n+p}$ to $X$ if and only if the following holds:
\beq &\forall 0 \leq i \neq j \leq n+p, \mbox{for any $b_i$ belonging to the label of some edge of $Q_i$},\\& \mbox{for any $b_j$ belonging to the label of some edge of $Q_j$, $b_i$ and $b_j$ commute}.\eeq

We can now give a proof of Theorem~\ref{thm:main1}, which we restate here.

\bthm $X$ is a locally CAT(0) cube complex, so $A(M)$ is cocompactly cubulated. \ethm

\bp Let $Q$,$Q'$, $Q''$ be cubes of $X$, which pairwise intersect in codimension $1$, and intersect globally in codimension $2$. Write $Q = \prod_{i=0}^{n+p} Q_i$, $Q' = \prod_{i=0}^{n+p} Q'_i$ and $Q'' = \prod_{i=0}^{n+p} Q''_i$.

Let $k,k',k'' \in \llb 0,n+p \rrb$ be such that
$$\forall i \neq k'', Q_i=Q'_i, \forall i \neq k', Q_i=Q''_i \mbox{ and } \forall i \neq k, Q'_i=Q''_i.$$

\bit
\item Assume first that $k=k'=k''$. Then the three cubes $Q_k$,$Q'_k$ and $Q''_k$ of $X_k$ pairwise intersect in codimension $1$ and globally intersect in codimension $2$. Since $X_k$ is locally CAT(0), there exists a cube $K_k$ in $X_k$ such that $Q_k$,$Q'_k$ and $Q''_k$ are codimension $1$ faces of $K_k$. Since for every $1 \leq i \leq n+p$, $X_i$ is a square complex and $K_k$ has dimension at least $3$, we deduce that $k=0$.

Let $K = K_0 \times \prod_{i=1}^{n+p} K_i$, where $\forall 1 \leq i \leq n+p, K_i=Q_i=Q'_i=Q''_i$.

\mk

We will check that the cube $K$ belongs to $X$: fix $0 \leq i \neq j \leq n+p$, and choose $b_i$ belonging to the label of some edge of $K_i$ and $b_j$ belonging to the label of some edge of $K_j$.
\bit
\item If $i,j \neq 0$, then $K_i=Q_i$ and $K_j=Q_j$, and since $Q$ is a cube of $X$, $b_i$ and $b_j$ commute.
\item If $i=0$ or $j=0$, assume that $i=0$. Then some edge of $K_0$ has label $\{b_0\}$. By definition of $X_0$, parallel edges in $K_0$ have the same labels, so $\{b_0\}$ is also the label of some edge of $Q_0$, $Q'_0$ or $Q''_0$: assume that $\{b_0\}$ is the label of some edge of $Q_0$. Since $Q$ is a cube of $X$, $b_0$ and $b_j$ commute.
\eit
As a consequence, $K$ is a cube of $X$.
\item Assume now that $k,k',k''$ are not all equal. Then $k,k',k''$ are pairwise distinct. Let $K= K_0 \times \prod_{i=1}^{n+p} K_i$, where for each $i \not\in \{k,k',k''\}$ we have $K_i=Q_i=Q'_i=Q''_i$, and also $K_k=Q'_k=Q''_k$, $K_{k'}=Q_{k'}=Q''_{k'}$ and $K_{k''}=Q_{k''}=Q'_{k''}$. Note that $Q_k$ has codimension $1$ in $Q'_k=K_k$, and similarly $Q'_{k'}$ has codimension $1$ in $K_{k'}$ and $Q''_{k''}$ has codimension $1$ in $K_{k''}$. Hence $K$ is a cube of the product $\prod_{i=0}^{n+p} X_i$, which contains each of $Q,Q',Q''$ with codimension $1$.

We will now check that $K$ is a cube of $X$. For any $0 \leq i \neq j \leq n+p$, assume that $b_i$ belongs to the label of some edge of $K_i$, and that $b_j$ belongs to the label of some edge of $K_j$. Without loss of generality, $i$ and $j$ are both different from $k$. So $K_i=Q_i$ and $K_j=Q_j$, and since $Q$ is a cube of $X$ we know that $b_i$ and $b_j$ commute. Hence $K$ is a cube of $X$.
\eit

According to Theorem~\ref{thm:gromov criterion}, $X$ is a locally CAT(0) cube complex.

\mk

The fundamental group of $X$ is given by its $2$-skeleton, and it is the quotient of the free product of $A(\Gamma|_{S_0}),\dots,A(\Gamma|_{S_{n+p}})$ obtained by adding the following commutation relations:
$$ \forall 0 \leq i \neq j \leq n+p, \mbox{if $a_i \in S_i$ and $a_j \in S_j$ commute in $A(M)$},\mbox{$a_i$ and $a_j$ commute in $\pi_1(X)$}.$$
The group $\pi_1(X)$ is therefore isomorphic to $A(M)$.

\mk

As a consequence, $A(M)$ is cocompactly cubulated.
\ep

We can now give the proof of Theorem~\ref{thm:main_dagger}, which we restate here.

\begin{thmE}
Conjecture~\ref{conj:main} holds for any Artin-Tits group satisfying property $(\dagger)$. In particular, Conjecture~\ref{conj:main} holds for Artin-Tits groups of type FC, and for Artin-Tits groups whose irreducible spherical parabolic subgroups have rank at most $2$.
\end{thmE}

\bp
Theorem~\ref{thm:main1} and Theorem~\ref{thm:main2}  precisely state that Conjecture~\ref{conj:main} holds for any Artin-Tits group satisfying property $(\dagger)$.

\mk

According to~\cite[Theorem~1]{godelle_cat0}, if the Deligne complex of an Artin-Tits group can be endowed with a piecewise Euclidean CAT(0) metric, then this Artin-Tits group satisfies property $(\dagger)$.

\mk

Consider an Artin-Tits group $A(M)$ of type FC. According to~\cite{charney_davis_kpi1}, the Deligne complex of $A(M)$, endowed with the cubical metric, is CAT(0). Therefore $A(M)$ satisfies property $(\dagger)$, and satisfies Conjecture~\ref{conj:main}.

\mk

Consider an Artin-Tits group $A(M)$ whose irreducible spherical parabolic subgroups have rank at most $2$. According to~\cite{charney_davis_kpi1}, the Deligne complex of $A(M)$, endowed with the Moussong metric, is CAT(0). Therefore $A(M)$ satisfies property $(\dagger)$, and satisfies Conjecture~\ref{conj:main}.
\ep

We can now prove Corollary~\ref{cor:braid}, which we restate here.

\begin{corG} The $n$-strand braid group $B_n$, or its central quotient $B_n/Z(B_n)$, is virtually cocompactly cubulated if and only if $n \leq 4$. \end{corG}

\bp We have seen that $B_3$ and $B_3/Z(B_3)$ are cocompactly cubulated.

\mk

Assume that $n \geq 4$, and let $\sigma_1,\dots,\sigma_{n-1}$ denote the standard generators of $B_n$. Since $m_{\sigma_1,\sigma_2}=m_{\sigma_2,\sigma_3}=3$, according to Theorem~\ref{thm:main2}, $B_n$ is not virtually cocompactly cubulated.

\mk

We will show that if the central quotient $B_n/Z(B_n)$ were virtually cocompactly cubulated, then $B_n$ itself would be virtually cocompactly cubulated. Assume that $G=B_n/Z(B_n)$ has a finite index subgroup $G_0$ acting geometrically on a CAT(0) cube complex $X$. Then $H=G_0Z(B_n)$ is a finite index subgroup of $B_n$ acting cocompatly on $X$. Consider the counting morphism $B_n \ra \Z$ sending each $\sigma_i$ to $1$. If we compose this morphism with the standard action of $\Z$ on $\R$ (seen as a cube complex) by translation, this defines a cocompact action of $H$ on $\R$. Since the centre $Z(B_n)$ acts properly on $\R$, we deduce that the action of $H$ on $X \times \R$ is proper and cocompact. So $B_n$ would be virtually cocompactly cubulated.

As a consequence, $B_n/Z(B_n)$ is not virtually cocompactly cubulated.
\ep

We now give a proof that Conjecture~\ref{conj:main2} implies Conjecture~\ref{conj:mainFC}.

\bp
Assume that an Artin-Tits group $A(M)$ is not of type FC, we want to prove that at least one of the three cases described in Conjecture~\ref{conj:main2} occur. Since $A(M)$ is not of type FC, there exist three pairwise distinct $a,b,c \in S$ such that $m_{ab},m_{bc},m_{ac}$ are finite and $\f1{m_{ab}} + \f1{m_{bc}} + \f1{m_{ac}} \leq 1$.

\mk

Assume first that at least one of $m_{ab},m_{bc},m_{ac}$ is odd, for instance $m_{ab}$ is odd. Then at least one of $m_{bc},m_{ac}$ is strictly bigger than $2$, for instance $m_{bc} >2$. Hence $a,b,c$ correspond to the first case described in Conjecture~\ref{conj:main2}.

Assume now that all of $m_{ab},m_{bc},m_{ac}$ are even. Then at most one of $m_{ab},m_{bc},m_{ac}$ is equal to $2$, for instance $m_{ab} \geq 4$ and $m_{ac} \geq 4$. So $a,b,c$ correspond to the second case described in Conjecture~\ref{conj:main2}.
\ep

To conclude, we give a proof of Corollary~\ref{cor:spherical}.

\bp
Assume first that $A$ is an Artin-Tits group of spherical type such that all its irreducible parabolic subgroups have rank at most $2$. Then $A$ is a direct product of cyclic groups and dihedral Artin groups, and according to Theorem~\ref{thm:main1} (or according to~\cite{brady_mccammond_artinthree}), $A$ is cocompactly cubulated.

\mk

Assume now that $A$ is an Artin-Tits group of spherical type containing an irreducible parabolic subgroup of rank at least $3$. So there exist three pairwise distinct $a,b,c \in S$ such that $m_{ab}=3$, $m_{ac} \neq \infty$ and $m_{bc} \neq 2$. The elements $a,b,c$ correspond to the first case in Conjecture~\ref{conj:main2}, so according to Theorem~\ref{thm:main_dagger}, $A$ is not virtually cocompactly cubulated.
\ep

\bibliographystyle{smfalpha_perso}
\bibliography{bibli}

\sign

\end{document}